\documentclass[11pt]{article}
\usepackage[margin=2cm]{geometry}
\usepackage[pagebackref=false,linktoc=page,colorlinks=true,linkcolor=blue,citecolor=blue]{hyperref}
\usepackage[T1]{fontenc}
\usepackage{natbib}
\usepackage{amsthm}
\usepackage{amssymb}
\usepackage{mathtools}
\usepackage{float}
\usepackage{caption}
\usepackage{color}
\usepackage{soul}
\usepackage{verbatim}
\usepackage{enumitem}
\usepackage{mathrsfs}
\usepackage{setspace}
\usepackage{pstool, psfrag} 
\usepackage[title]{appendix}
\usepackage[retainorgcmds]{IEEEtrantools}
\usepackage{bm}
\usepackage{amsmath}
\usepackage{graphicx}
\usepackage{authblk}
\renewcommand{\&}{and}
\newtheorem{theorem}{Theorem}
\newtheorem{lemma}{Lemma}

\newtheorem{remark}{Remark}
\newtheorem*{theorem*}{Theorem}
\newtheorem{proposition}{Proposition}

\newcommand{\meas} {\ensuremath{\nu}}

\newcommand{\GG} {\ensuremath{\mathcal{G}}}
\newcommand{\cind}{\ensuremath{\,\overset {d }{\longrightarrow }\,}}
\newcommand{\wk}{\ensuremath{\,\overset {\mathrm {w} }{\longrightarrow
    }\,}}

\newcommand{\dd}{\ensuremath{\mathrm{d}}}

\newcommand{\AAA}{\ensuremath{A}_0}
\renewcommand{\AA}{\ensuremath{A}_1}
\newcommand{\AB}{\ensuremath{A}_2}

\newcommand{\BA}{\ensuremath{B}_1}

\newcommand{\HR}{H\"{u}sler-Reiss}

\newcommand{\PP}{\ensuremath{\mathbb{P}}}
\newcommand{\RR}{\ensuremath{\mathbb{R}}}
\newcommand{\NN}{\ensuremath{\mathbb{N}}}
\newcommand{\SSS}{\ensuremath{\mathbb{S}}}
\newcommand{\EE}{\ensuremath{\mathbb{E}}}
\newcommand{\level}{\ensuremath{t}}
\newcommand{\thres}{\ensuremath{t}}
\newcommand{\ttime}{\ensuremath{i}}
\newcommand{\sm}{\ensuremath{\setminus}}
\newcommand{\Vbnd}{\ensuremath{N_{\GG}}}

\newcommand{\Gbnd}{\ensuremath{\mathcal{N}_{\GG}}}
\newcommand{\condv}{\ensuremath{v}}
\newcommand{\dims}{\ensuremath{s}}
\newcommand{\expdist}{\ensuremath{F^1_E}}

\newcommand{\bigCI}{\mathrel{\text{\scalebox{1.07}{$\perp\mkern-10mu\perp$}}}}

\providecommand{\keywords}[1]
{  \textbf{\text{Keywords}} #1
}
\title{\textbf{Decomposable tail graphical models}}
\author[$\dagger$]{\textbf{Adrian Casey} }
\author[$\dagger$\footnote{Corresponding author}]{\textbf{Ioannis
    Papastathopoulos}}
\affil[$\dagger$]{\small School of Mathematics and Maxwell Institute,
  University of Edinburgh, Edinburgh, EH9 3FD}
\affil[$$]{\small
  a.casey2@exseed.ed.ac.uk 
  $\quad$ i.papastathopoulos@ed.ac.uk} 

\date{}

\begin{document}
%
%
\setstretch{1.15}
\maketitle
\begin{abstract}
  We develop an asymptotic theory for extremes in decomposable
  graphical models by presenting results applicable to a range of
  extremal dependence types.\ Specifically, we investigate the weak
  limit of the distribution of suitably normalised random vectors,
  conditioning on an extreme component, where the conditional
  independence relationships of the random vector are described by a
  chordal graph.\ Under mild assumptions, the random vector
  corresponding to the distribution in the weak limit, termed the
  \textit{tail graphical model}, inherits the graphical structure of
  the original chordal graph.\ Our theory is applicable to a wide
  range of decomposable graphical models including asymptotically
  dependent and asymptotically independent graphical models.\
  Additionally, we analyze combinations of copula classes with
  differing extremal dependence in cases where a normalization in
  terms of the conditioning variable is not guaranteed by our
  assumptions.\ We show that, in a block graph, the distribution of
  the random vector normalized in terms of the random variables
  associated with the separators converges weakly to a distribution we
  term \textit{tail noise}.\ In particular, we investigate the limit
  of the normalized random vectors where the clique distributions
  belong to two widely used copula classes, the Gaussian copula and
  the max-stable copula.
\end{abstract}

\section{Introduction}
Assessing risk from extreme multivariate events is pivotal in a broad
range of applications ranging from climate science and hydrology to
finance and insurance.\ The dependence structure of extremes, however,
can be highly intricate and their probabilistic and statistical
modelling is a major challenge.\ Even in two dimensions, modelling
extremes can be arduous since there are two possibilities for the
types of extremal association between two random variables.\ On one
hand, \textit{asymptotic independence}, through which two variables
may be strongly dependent at moderately extreme levels, but
ultimately, the largest extremes occur independently; or
\textit{asymptotic dependence}, where the largest extremes may occur
jointly.\ Such challenges are exacerbated when considering extremes of
real-valued random vectors $\bm X=(X_1,\dots,X_d)$ in high dimensions
$d\gg 2$ since, not only do the types of possible extremal dependence
between random variables explode combinatorially, but also because
data in joint tail regions are scarce.\ The information that is
typically available in joint tail regions renders unreliable both
empirical estimates and statistical models with saturated or
unstructured dependence forms.\ We address these challenges by
focusing on the class of decomposable graphical models, a class of
parsimonious and flexible statistical models that can scale with the
dimension $d$ of the random vector $\bm X$.\ These models also allow
for clique distributions exhibiting asymptotic dependence or
asymptotic independence, also combined within the same graphical
model.

Studying multivariate extreme values on simple graphical structures
dates back to the work of \cite{smit92} who investigated the scaling
properties of positive-recurrent real-valued Markov chains, subject to
the assumption that the marginal distributions of adjacent random
variables in the Markov chain belong to the domain of attraction of a
bivariate max-stable distribution.\ \cite{smit92} showed that when the
Markov process is conditioned to exceed a high value in its initial
state, then the rescaled Markov process converges weakly, as the level
of extremity tends to the upper end point of the marginal distribution
to a random walk, termed the \textit{tail chain}.\ The limit
process inherits the graphical structure from the original Markov
chain and this has the useful practical implication of justifying
structured limit processes as candidate statistical extreme value
models.\ Similar asymptotic theory for more general $\RR^d$-valued
Markov processes with extreme initial states was given by
\cite{perf97}, \cite{yun1998extremal} and \cite{janssege14}.\ 

The subject of extremes in graphical structures has seen a rich vein
of research recently.\ \cite{engelke2020graphical} introduced a
general theory for conditional independence for multivariate
generalized Pareto distributions and opened up a way for formulating
statistical extreme value models that maintain sparsity, termed
extremal graphical models.\ \cite{engelke2020graphical} and
\cite{engelke2022structure} demonstrate methods of inference and
graphical learning for extremal graphical models.\ In these papers as
in \cite{segers2020one, asenova2020inference,
  engelke2021sparse,engelke2022structure} and
\cite{hentschel2022statistical} there is a concentration on tractable
models using the \HR \ generalized Pareto distribution.\ Multivariate
generalized Pareto distributions arise as the only non-trivial limit
distributions of renormalized multivariate exceedances
\citep{roota:06} and have a number of attractive theoretical and
practical properties.\ Firstly, they are the multivariate analogue of
the generalized Pareto distribution \citep{pick75} in the sense that
they are the only threshold-stable multivariate distributions
\citep{kiriliouk2019peaks}.\ If a random vector follows a multivariate
generalized Pareto distribution, then the conditional distribution of
a multivariate exceedance is a location-scale transformation of the
original distribution.\ The practical relevance of the
threshold-stability property is that the form of the model does not
change at higher levels and this is useful when extrapolating into the
tail.\ Secondly, statistical models based on multivariate generalized
Pareto distributions often lead to simple and easily computable
likelihood functions, which is in stark contrast to their associated
multivariate max-stable distributions, which arise as the only
non-trivial limit distributions of renormalized componentwise maxima
of i.i.d. random vectors.

A crucial assumption of all of the aforementioned approaches, however,
is that of the graphical model under consideration being multivariate
regular varying in the interior of the cone $\RR^d_+\sm\{\bm 0\}$
\citep{resn07}, which implies that all such processes exhibit
asymptotic dependence at all lags in a Markov chain, or between all
random variables in a graphical model.\ Although powerful for
modelling extremes of univariate random variables, this assumption,
which holds for the multivariate peaks-over-threshold analogue and its
counterpart approach based on componentwise maxima presupposes a
fundamental \textit{dependence assumption}, that is, both approaches
assume that $\chi_V>0$ where the coefficient of asymptotic dependence
$\chi_A$ is more generally defined for any non-empty
subset$A \subseteq V$ of random variables by
\begin{equation} 
  \chi_{A} = \lim_{q\rightarrow 1} \, \PP(F_i(X_i) > 1-q\,:\,i\in
  A)/(1-q),
  \label{eq:chi}
\end{equation}
where $F_i$ denotes the distribution function of $X_i$.\ This
dependence assumption implies that all variables may attain their
largest values simultaneously, regardless of how large the dimension
$d$ is.\ The frequent absence of this form of extremal dependence in
observed data motivates us to study asymptotic theory for a wider
variety of models.

We focus on graphical models with respect to conditional independence
graphs that are connected and chordal
Our approach rests on conditioning on an arbitrary variable
$X_{\condv}$ exceeding a high threshold and studying the limiting
behaviour of the random vector, affinely renormalized with functions
of the conditioning variable.\ The key difference and strength of our
approach is on the choice of exponentially-tailed marginals, which
allows to study affine renormalizations and uncover structure that can
be missed by rescaling in heavy-tailed margins.\ Our approach is based
on the asymptotic theory of Heffernan--Resnick--Tawn
\citep{hefftawn04, heffernan2007limit} and we follow a similar program
to that in \cite{papastathopoulos_strokorb_tawn_butler_2017}, where a
unified theory of both fully asymptotically dependent and fully
asymptotically independent time-homogeneous real-valued Markov chains
was developed, and where it was shown, assuming exponentially tailed
margins and using affine renormalization, that a scaled autoregressive
process is obtained as the weak limit of the renormalized Markov
chain, a result which generalizes the concept of the tail chain to the
hidden tail chain.\ In what follows, we use the generic term
\textit{tail graphical model} to denote the random vector in the
normalized limit, but mention here that our work is essentially an
analogue of the hidden tail chainthat was studied by
\cite{papastathopoulos2023hidden}.\ Under mild assumptions, we show
that the weak limit exists for the affinely normalized random vector
and that the limit random vector inherits the graphical structure of
the random vector $\bm X$.\ We study widely used examples of
decomposable graphical models and explain this structure in detail for
distributions with max-stable H\"{u}sler-Reiss and Gaussian random
subvectors associated with each maximal clique, which are
asymptotically dependent and asymptotically independent, respectively.

In this work, we also develop a similar approach for graphical models
with clique distributions that exhibit different types of extremal
dependence.\ For the sake of simplicity, we restrict the analysis to
block graphs, that is, decomposable graphs where maximal cliques
intersect at a single vertex.\ For such cases, we identify situations
where standard asymptotic theory breaks and consequently, we introduce
the concept of \textit{tail noise}, which describes the weak limit
that emerges when the random vector $\bm X$ is normalized with
functions of the variable associated with the intersection rather than
with functions of the conditioning variable.\ To the best of our
knowledge, this is the first work that considers both extremal
dependence types in a graphical models and hence, it extends previous
work by encompassing a wider range of decomposable graphical models.

The paper is structured as follows. In Section \ref{sec:background} we
introduce background material on graphical models and extreme value
theory.\ In particular, Section \ref{sec:MRV} discusses multivariate
regular variation and presents two preliminary findings, one that
identifies Markov random fields with multivariate regularly varying
clique distributions as elements in the domain of attraction of
extremal graphical models, and the other which provides conditions
subject to which asymptotically dependent cliques may join with
asymptotically independent cliques on a conditional independence
graph.\ In Section \ref{sec:ce} we describe conditional extreme value
theory which we use subsequently in Section \ref{sec:main}. Section
\ref{sec:main} presents our main results theorems, one associated with
the weak convergence of affinely renormalized graphical models to a
tail graphical model and the other with Section \ref{sec:examples}
contains examples of the analysis applied to two widely used
multivariate distributions.\ Relevant material on graph theory is
summarized in Appendix \ref{app:graph}.  All proofs are presented in
Appendices~\ref{app:proofs} and \ref{app:exproofs}.

\subsection{Some notation}
The set $V$ denotes the sequence of integers $\{1, \dots, d \}$.\ The
notation $|A|$ denotes the cardinality of the set $A$, for example,
$\lvert V\rvert = d$.\ The notation $ C_b(\RR^d)$ denotes the space of
continuous and bounded real functions on $\RR^d$.\ If measures
$\mu_\level$ and $\mu $ on $\RR^d$ satisfy
$\mu_\level f\rightarrow \mu f$ as $\level\rightarrow \infty$ for
every $f\in C_b(\RR^d)$, we say $\mu_\level\wk \mu$.\ We write
$\bm X_\level \cind \bm X$ when $F_\level\wk F$, where $F_\level$ and
$F$ denote the distributions of $\bm X_\level$ and $\bm X$.\ We use
the notation \smash{$\mathcal{N}_{d}\big ( \bm \mu, \bm\Sigma \big )$
} to denote the $d$ dimensional multivariate normal distribution with
mean $\bm \mu$ and covariance $\bm\Sigma $. For an arbitrary non-empty
subset $A\subset V$ we let $\bm x_A=(x_v\,:\,v \in A)$, so
$\bm x_A \in \RR^{|A|}$.\ An $|A|\times |A|$ by matrix is denoted as
any of $\bm \Sigma=\bm \Sigma_{A}=\bm \Sigma_{AA}$ whereas for two
arbitrary non-empty subsets $D$ and $E$ of $V$, we let
$\bm \Sigma_{D,E}=(\sigma_{v v'})_{v\in D, v'\in E}$ denote a
$|D|\times|E|$ submatrix of $\bm \Sigma$.\ All arithmetical operations
on vectors are interpreted as element-wise and by convention, if
$f \,:\,\RR\rightarrow \RR$ and $\bm x\in\RR^d$, $d\in\NN$, then
$f (\bm x) = (f (x_1), \ldots, f (x_d))$.\ If $\bm x$ is a $d$
dimensional vector indexed by the vertices $V$ of a graph $\GG$ then
$\bm x_{V \sm v}$ denotes the $d-1$ dimensional vector indexed by the
subgraph induced by the vertices $V \sm v$.\ Similarly, for
$S\subset V$, $\bm x_{V \sm S}$ denotes the $d-|V \cap S|$ dimensional
vector indexed by $V\sm S$.\ We shall extend this notation to include
such technically inappropriate expressions as $\bm X_{C \sm \condv}$
and $\bm X_{C \sm S}$ where the indexing is over the vertices of the
subgraphs $C \sm \condv$ and $C \sm S$ respectively. A measurable
function $g\,:\,\RR_+\rightarrow \RR_+$ is regularly varying at
infinity with index $\rho \in \RR$, denoted by $g\in\text{RV}_{\rho}$
if, for any $x>0$, $g(\level x)/g(\level)\rightarrow x^{\rho}$ as
$\level\rightarrow \infty$.

\section{Background}
\label{sec:background}
\subsection{Decomposable graphical models}
We summarize here the most important properties of decomposable
graphical models that we shall use throughout.\ Further introductory
material and terminology on graph theory is provided in Appendix
\ref{app:graph}.\ A random vector $\bm X=(X_1,\dots,X_d)$ is said to
follow a graphical model with conditional independence graph
$\GG=(V,E)$, if its distribution $F_{\bm X}$ satisfies the pairwise
Markov property relative to $\mathcal{G}$, that is,
$X_i\bigCI X_j \mid \bm X_{V\setminus \{i,j\}}$ for all
$(i,j)\notin E$.\ A graphical model is termed decomposable if its
conditional independence graph $\GG$ is chordal.\ For a connected
chordal graph, we can always order the set of (maximal) cliques as
$\{{C}_1, \dots, C_N\}$ such that, for all $i=2, \dots, N$,
\begin{equation} S_\ttime := C_\ttime \cap ( \cup_{j=1}^{i-1} C_j ) \subset C_k
  \quad \text{for some } \, k < \ttime.
  \label{eq:ri}
\end{equation}
This condition is called the running intersection property and from
now on, we shall always work with connected chordal graphs, for which
the set of cliques has been ordered so that expression \eqref{eq:ri}
holds true. We denote this ordered set by $\mathcal{M}(\GG)$.\

\begin{figure}[htbp!]
  \includegraphics[scale=0.55,trim=50 100 50 30]{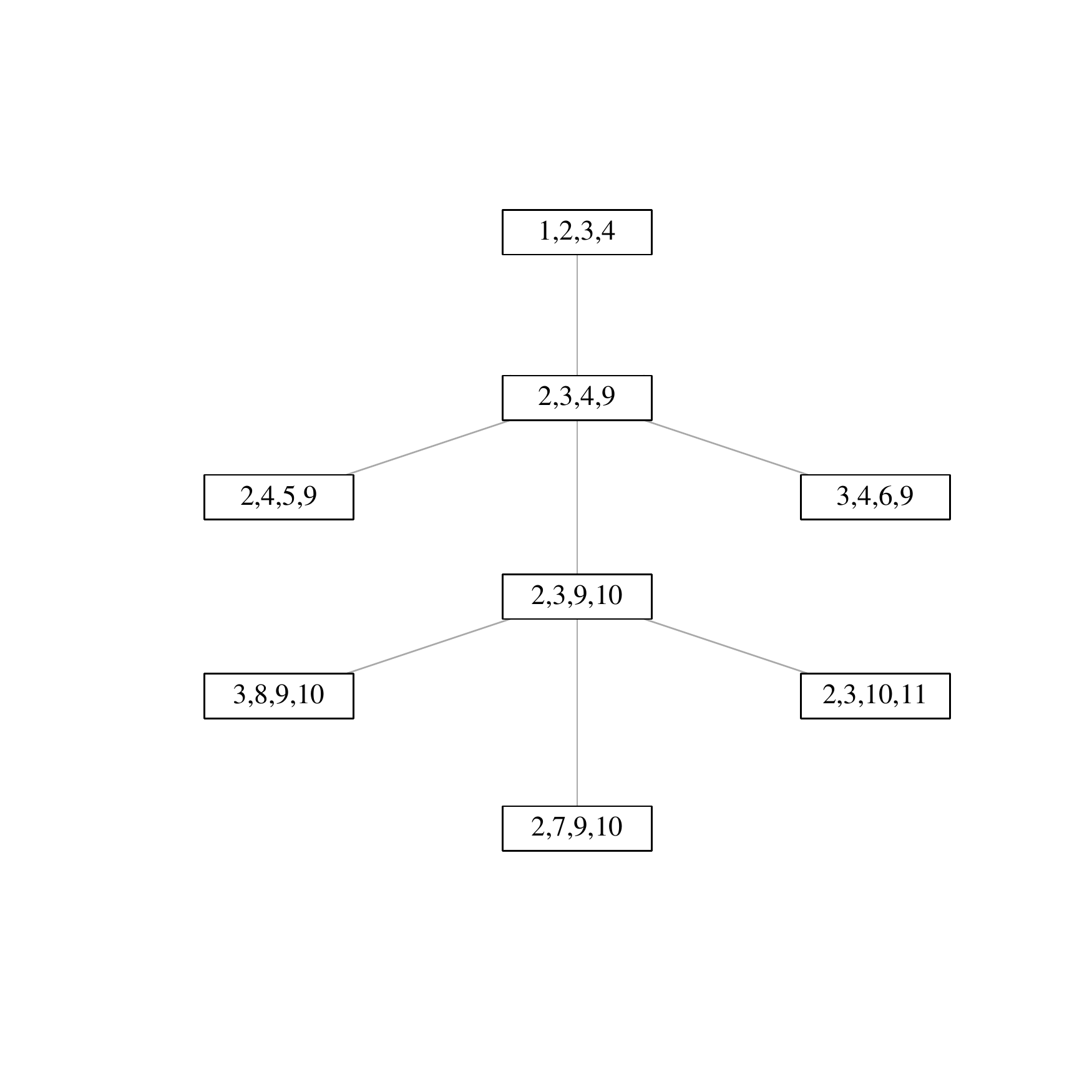}
  \includegraphics[scale=0.55, trim=50 100 50 30]{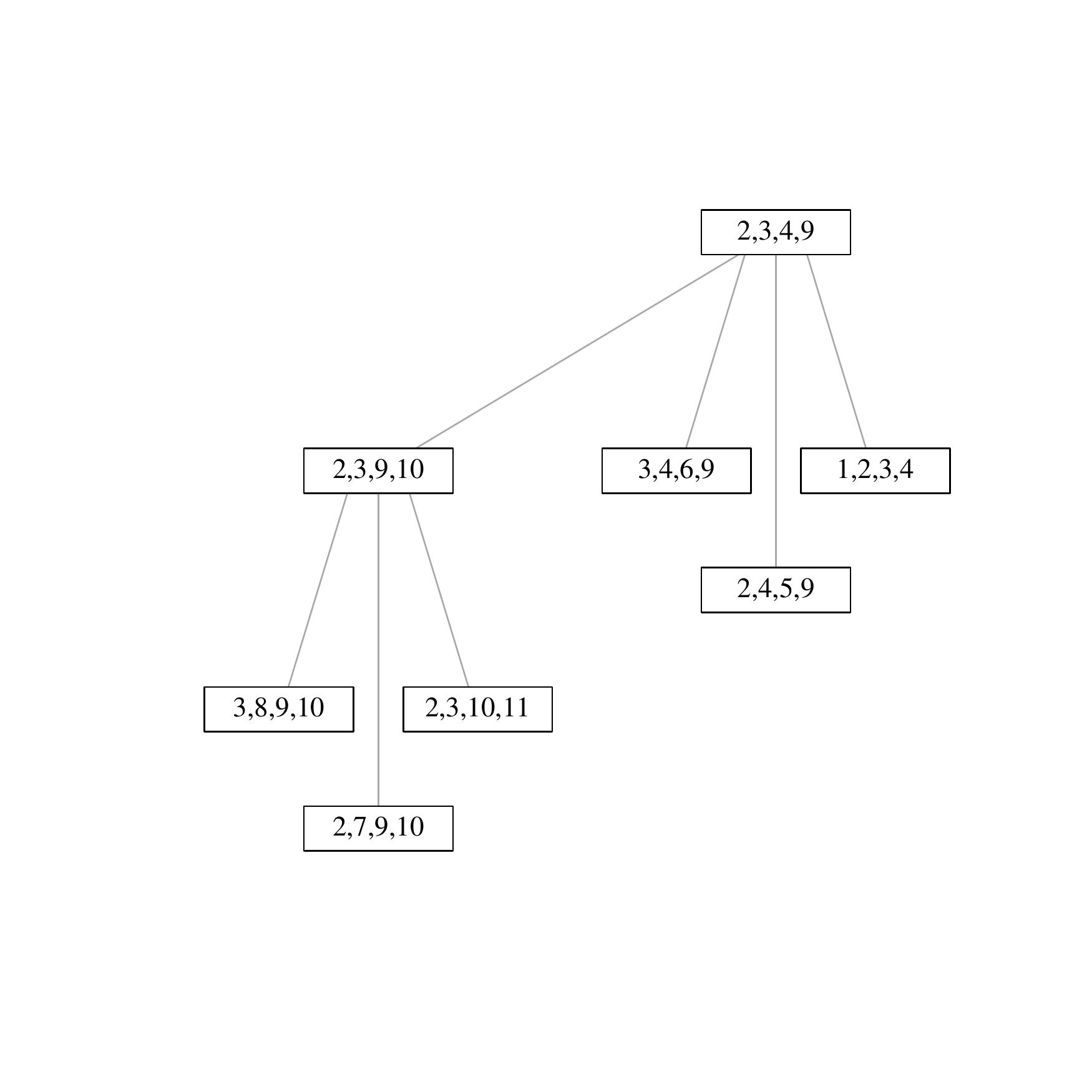}
  \caption{The junction tree of the Goldner-Harary graph (Figure 3)
    with two different cliques chosen as the root.  The root is shown
    at the top.}
  \label{fig:arbroot}
\end{figure}

\begin{figure}[htbp!]
  \includegraphics[scale=0.4, trim=180 150 180 60]{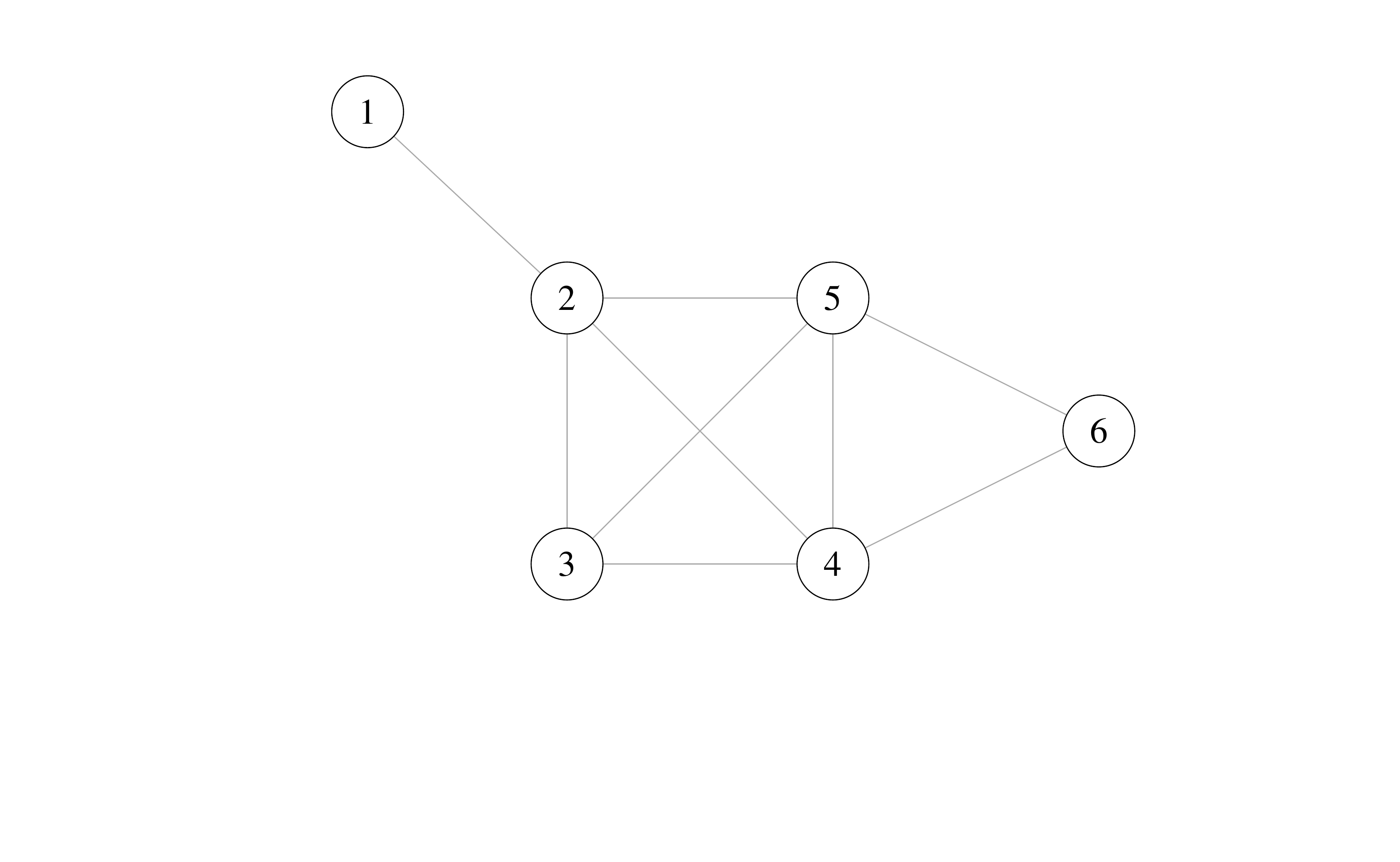}
  \includegraphics[scale=0.4, trim=180 150 180 60]{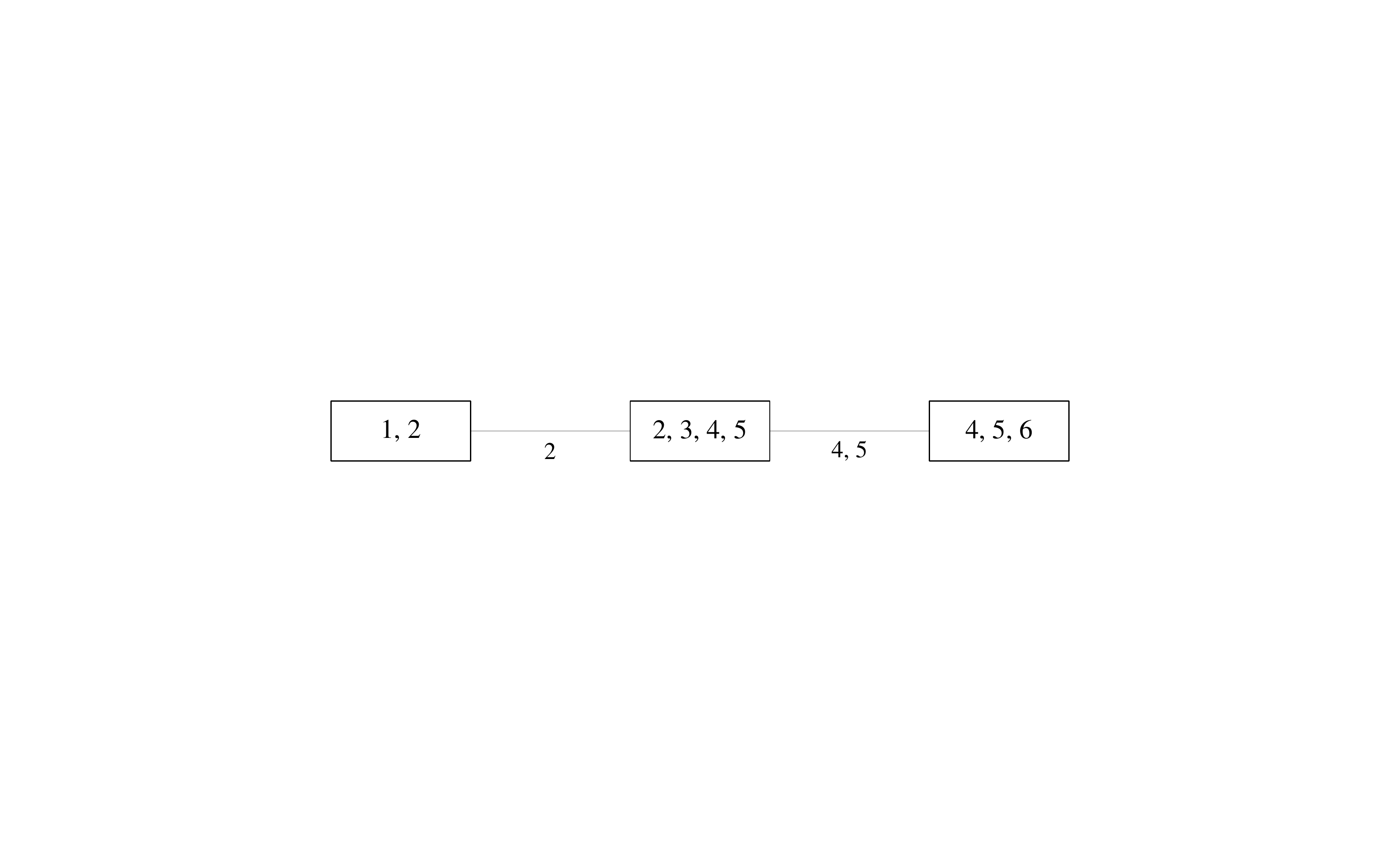}
  \caption{\textit{Left:} Example of a decomposable graph.\
    \textit{Right:} its junction tree.\ The vertices in the junction
    tree are the cliques of the graph and the labels on the edges
    display the separators.}
  \label{fig:sepsimple}
\end{figure}

\begin{figure}[htbp!]
  \includegraphics[scale=0.4, trim=200 0 160 20]{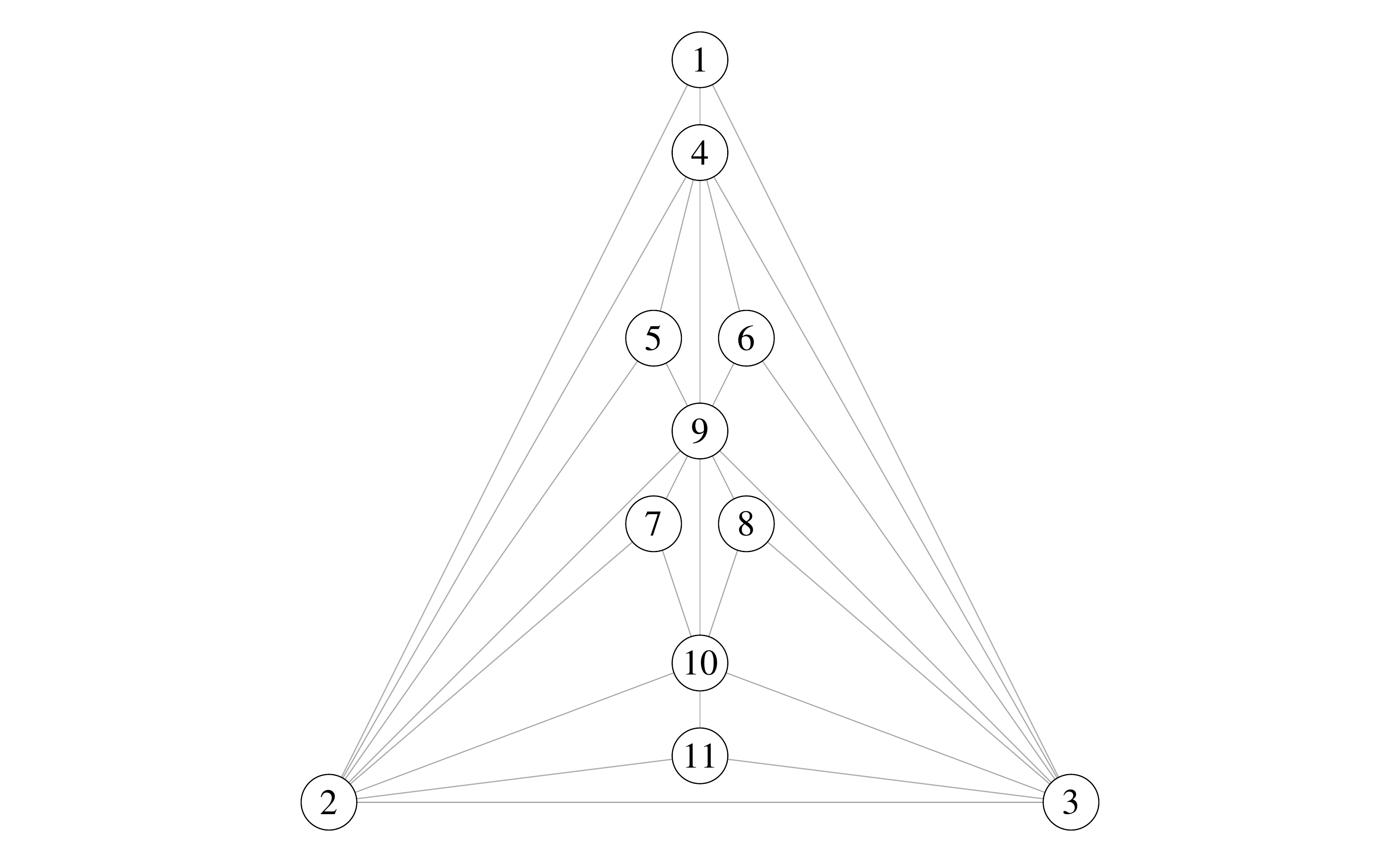}
  \includegraphics[scale=0.4, trim=200 0 160 20]{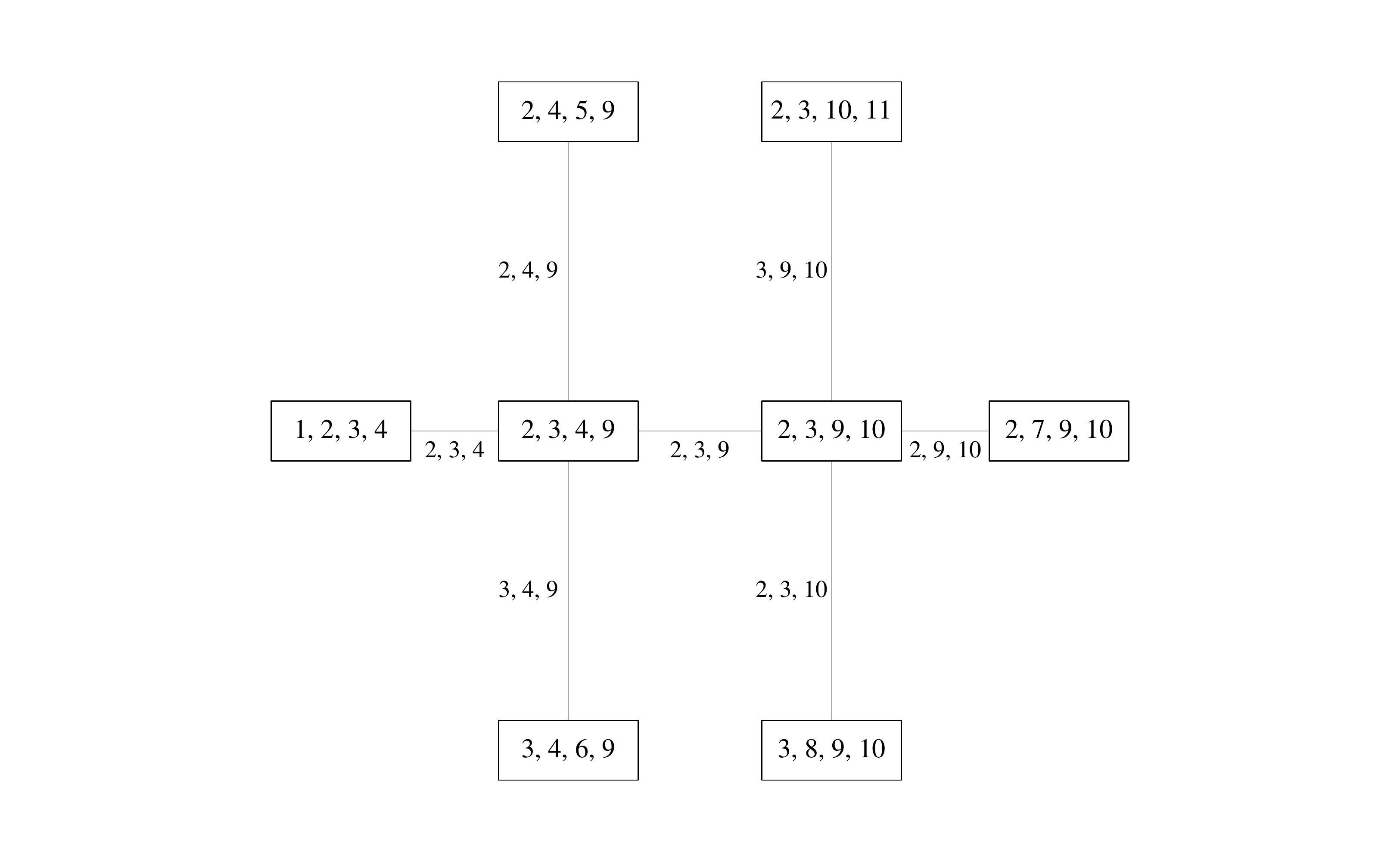}
  \caption{\textit{Left}: The Goldner-Harary graph.\ \textit{Right:}
    its junction tree.\ The nodes in the junction tree are the
    cliques of the graph and the labels on the edges display the
    separators.}
  \label{fig:sepgh}
\end{figure}

It is useful to regard the set $\mathcal{M}(\GG)$ as vertices of a
tree, which we term nodes, and with an edge between nodes if
$C_i \cap C_j \neq \varnothing$ for any
$C_i, C_j \in \mathcal{M}(\GG)$.  Such a tree, with $\mathcal{M}(\GG)$
as nodes, is called a junction tree, denoted by $\GG^*$, if for any
two distinct nodes $C_i$ and $C_j$, every node on the unique path
connecting $C_i$ and $C_j$ in the tree represents a clique that
contains the intersection $C_i \cap C_j$.\ A graph $\GG$ is chordal if
and only if there exists a junction tree $\GG^*$ \citep[see
Proposition 2.27]{laur96}.\ Two examples of chordal graphs and their
junction trees are shown in Figures \ref{fig:sepsimple} and
\ref{fig:sepgh}.
\begin{remark}\normalfont
  The choice of $C_1$ in the running intersection property is
  arbitrary \citep[see Lemma 2.18]{laur96}.  Equivalently, in
  representing the cliques on a junction tree, we can select an
  arbitrary node as a root see for example Figure \ref{fig:arbroot}.
  \label{rem:rip}
\end{remark}
An important restriction on decomposable graphical models is the
requirement of compatibility for the distributions of subvectors, in
the sense that for any two cliques $C, C' \in \GG^*$ with
$S = C\cap C' \neq \emptyset$, then the distributions $F_{\bm X_{C}}$
and $F_{\bm X_{C'}}$, when marginalized over the variables indexed by
$C \sm S$ and $C \sm S'$, must each yield the same distribution
$F_{\bm X_{S}}$.\ When such a restriction holds then there is a unique
distribution $F_{\bm X_{C \cup C'}}$ \citep{dawid1993hyper}.\ If a
random vector has conditional independence properties described by a
graph with singleton separators, termed a block graph, then this
compatibility implies that the univariate marginal distributions
associated with the singleton separator sets are identical.
If the graphical model $\bm X$ has a strictly positive joint density
$f_{\bm X}(\bm x)$, where $\bm x \in \RR$ then the density factorizes
according to
\begin{IEEEeqnarray*}{rCl}
  f_{\bm X}(\bm{x}) &=& \prod_{i=1 }^N f_{\bm
    X_{C_i}}(\bm{x}_{C_i})\Big / {\prod_{i=2 }^N f_{\bm
      X_{S_i}}(\bm{x}_{S_i})}, \quad \bm x \in \RR^d,
  \label{eqn:dense0}
\end{IEEEeqnarray*}
where $f_{\bm X_A}$ denotes the density of the random vector
$\bm X_A$, with $A \subseteq V$, $A \neq \emptyset$ \citep{laur96}.\
We use this result in Section \ref{sec:MRV} in the context of models
with regularly varying clique distributions that admit a Lebesgue
density.  In our more general results we do not assume the existence
of a joint density, but instead rely only on the existence of Markov
transition probability kernels of form
$\PP(\bm X_{C_i \sm S_i} \in A \mid \bm X_{S_i} = \bm x_{S_i})$ where
$C_i$ is a clique and $S_i$ a separator joining $C_i$ to a
neighbouring clique, and that,
\begin{equation*}
  \PP(\bm X_{C_i \sm S_i} \in A \mid \bm X_{V \sm \{C_i \sm S_i\}} = \bm x_{V \sm \{C_i \sm S_i\}}) = \PP(\bm X_{C_i \sm S_i} \in A \mid \bm X_{S_i} = \bm x_{S_i}) \quad i=2, \dots, N.
\end{equation*}
By construction, the existence of these transition probability kernels
ensures the conditional independence relationships implied by the
graph $\GG= ( V , E)$.
Decomposable graphs provide a wide and varied subset of all graphical
models, and can encode a rich dependency structure.\ We illustrate
this with two graphs, one is a simple example, (Figure
\ref{fig:sepsimple}), and the other, the Goldner-Harary graph (Figure
\ref{fig:sepgh}), shows the variety of dependence structures that can
be described by a decomposable graph.\ The Goldner-Harary graph and
its junction tree also demonstrate the power of the junction tree as
an analytical tool for illustrating the conditional independence
between cliques in a clear and simple form.

\subsection{Multivariate regular variation}
\label{sec:MRV}
\color{black} Multivariate regular variation is the most frequently used
assumption for the probabilistic behaviour of the distributional tails
of a random vector and is summoned in a wide variety of statistical
approaches for multivariate extremes.\ \color{black}A $d$-dimensional
random vector $\bm X $ is multivariate regularly varying on the cone
$\EE=[0,\infty)^{V}\sm\{\bm 0_{|V|}\}$, with index $\alpha>0$, if for
any relatively compact $B\subset \EE$, there exists a limit measure
$\nu$ such that
\begin{equation}
  \level \,\PP\left(\bm X / b(\level) \in B\right) \rightarrow \meas(B),
  \quad \text{as $\level\rightarrow\infty$},
  \label{eq:MRV}
\end{equation}
with $\meas(\partial B) = 0$, $b(\level)\in \text{RV}_{1/\alpha}$. The
limit measure $\meas$ is necessarily $(-\alpha)$-homogeneous.\
Convergence \eqref{eq:MRV} is most conveniently expressed in
standardized Fr\'{e}chet$(1)$ margins which we assume in what follows,
that is, we assume that the marginal distributions of the random
vector $\bm X$ satisfy $F_{X_v}(x)=\exp(-1/x)$, $x>0$, for all
$v \in V$.\ This assumption about the form of the margins entails that
$b(\level)=\level$ and $\alpha=1$.\ The set $B$ is often taken as
$[\bm 0_{|V|},\bm x]^c$, $\bm x\in \RR_+^{V}$, which leads to the exponent
measure defined by $\Lambda(\bm x)=\meas([\bm 0_{|V|}, \bm x]^c)$.\ Subject
to this choice for the margins and due to the homogeneity of the limit
measure $\nu$, multivariate regular variation \eqref{eq:MRV} is
usually expressed in terms of the radial-angular decomposition
\[
  \PP\left(\frac{\bm X}{\lVert \bm X\rVert}\in A, \lVert\bm X\rVert >
    \level\, r ~\Big|~ \lVert\bm X\rVert > \level\right) \rightarrow
  H(A)\, r^{-1}, \quad r\geq 1,
\]
where $\lVert\,\cdot\,\rVert \,:\,\RR^{|V|}\rightarrow\RR_+$ is any norm,
$A\subset \SSS_{|V|-1}=\{\bm w\in\RR_+^{|V|}\,:\,\lVert \bm w \rVert=1\}$
and $H$ is a Radon measure on $\SSS_{|V|-1}$, termed the
spectral measure of the limit exponent measure $\Lambda$.\ %

Multivariate regular variation implies domain of attraction properties
for the distribution $F_{\bm X}$.\ Let
$\{(X_v^i\,:\,v\in V)\,:\,i=1,\dots,n\}$ denote a random sample from
$F_{\bm X}$ and let
$\bm M_n = \left ( \bigvee_{i=1}^n X_v^{i}\,:\,v\in V \right)$ denote
the vector of component-wise maxima.\ Then, $F_{\bm X}$ belongs to the
domain of attraction of a max stable distribution, that is,
\begin{equation}
 \lim_{n \rightarrow \infty} \text{Pr} \left(\bm M_n \leq n\,\bm y
  \right ) =\exp \left \{ - \Lambda (\bm y) \right \},
  \label{eq:expfn}
\end{equation} 
for $\bm y\in \EE$ which are continuity points of the limit.\ Here and
below, $\bm x \leq \bm y$ means that $\bm y - \bm x \in \EE$ and
$\bm x < \bm y$ means that $\bm x \leq \bm y$ and $\bm x\neq \bm y$.\
Similarly, the vector of normalized multivariate exceedances of the
random vector $\bm X$ \citep{rootzen2018multivariate} converges in
distribution to a multivariate Pareto distribution, that is,
\begin{equation}
  \lim_{\level\rightarrow \infty} \text{Pr} \left ( \bm X \leq \level\, \bm y \ ~\Big |~ \lVert \bm X \rVert_{\infty} > \level \right ) = \frac{\Lambda \left(\bm y\wedge \bm 1   \right) - \Lambda \left (  \bm y   \right) } { \Lambda \left ( \bm 1  \right )},
  \label{eq:pareto}
\end{equation}
for $\bm y\in \EE$ which are continuity points of the limit.

The parts of $\EE$ where $\nu$ places mass, or equivalently, the
faces of $\SSS_{|V|-1}$ where $H$ places mass, illustrate the
broad-scale extremal dependence structure of $\bm X$.\ The cone $\EE$
can be decomposed into the disjoint union
$\EE=\bigcup_A\EE_{A\subseteq V, A\neq \emptyset}$, where
\[
\EE_A = \{\bm x \in \EE\,:\, \bm x_A > \bm 0_{|A|}\quad \text{and}\quad \bm x_{V\sm A}=\bm 0_{|V \sm A|}\}.
\]
If $\nu(\EE_A) > 0$, then the variables indexed by $\bm X_A$ can take
their most extreme values simultaneously, while the components of
$X_{V\sm A}$ are non-extreme.\ In what follows, we adopt the same
notation for denoting the parts of the lower dimensional cones
$\mathbb{E}^M:=[0,\infty)^{M}\sm\{\bm 0\}$, where $M\subset V$ with
$M\neq \emptyset$, so that, for any $A\subseteq M$,
$\mathbb{E}_A^M=\{\bm x \in \EE^M\,:\, \bm x_A > \bm 0_{|A|}\quad
\text{and}\quad \bm x_{M\sm A}=\bm 0_{|M \sm A|}\}$

A common assumption that underpins a wide range of existing approaches
for characterizing extremes of graphical models \citep{kulik2015heavy,
  engelke2020graphical, segers2020one, asenova2021extremes}, is to
assume that the random vector $\bm X$ is multivariate regularly
varying with a limit measure $\nu$ that places all of its mass in the
interior $\EE_V$ of $\EE$, in the sense that $\nu(\EE\sm\EE_V)=0$ (or
equivalently, with a spectral measure that satisfies
$H(\SSS_{|V|-1}\sm \partial \ \SSS_{|V|-1})=0$).\ With this
assumption, and an additional requirement that the clique
distributions of $\bm X$ have Lebesgue densities, we can identify the
density of the exponent measure $\Lambda$ that represents the
max-stable domain of attraction of $\bm X$.
\begin{proposition}
  \label{prop:maxstable}
  Let $\bm X$ be a graphical model on a decomposable graph
  $\GG$ admitting a Lebesgue density $f_{\bm X}$ on $\RR_+^d$.\
  Suppose that for every clique $C_i$ of $\GG$, $\bm X_{C_i}$
  is multivariate regularly varying on
  $\EE^{C_i}=[0,\infty]^{|C_i|}\sm\{\bm 0\}$ with limit measure
  $\nu^{C_{i}}$ satisfying
  \begin{enumerate}
  \item[$(i)$]$ \nu^{(C_i)}(\EE^{C_i}\sm\EE_{C_i}^{C_i})=0$;
  \item[$(ii)$] The exponent measure
    $\Lambda^{(C_i)}(\bm y)=\nu^{(C_i)}([\bm 0, \bm x]^c)$ is absolutely
    continuous with respect to the Lebesgue measure on $\EE^{C_i}$ with
    density $\lambda^{(C_i)}$, in the sense that
    $\Lambda^{(C_i)}(\bm y) = \int \lambda(\bm x_{C_i})\,I(\bm
    x_{C_i}\in{[\bm 0,\bm y_{C_i}]^c}) ~\dd\bm x_{C_i}$.
  \end{enumerate} 
  Assume further that as $n\rightarrow\infty$
  \begin{equation}
    n^{|C_{\ttime}|+1}f_{\bm X_{C_i}}(n \, \bm x_{C_i})\rightarrow    
    \lambda^{(C_i)}(\bm x_{C_{i}}), \quad \text{$1,\dots, N$},
    \label{eq:diff_doa}
  \end{equation}
  outside a set of Lebesgue measure zero.\ Then $F_{\bm X}$ belongs to
  the domain of attraction of a max-stable distribution (equivalently
  $F_{\bm X}$ belongs to the domain of attraction of a multivariate
  Pareto distribution).\ The exponent measure $\Lambda$ of the limit
  multivariate max-stable and multivariate Pareto distribution is
  absolutely continuous with respect to Lebesgue measure on $\EE$ and
  has density 

  \begin{IEEEeqnarray}{rCl}
    \lambda(\bm y)&=&  \lambda^{(C_1)}(\bm y_{(C_1)}) \,\prod_{i= 2}^N
    \{\lambda^{(C_i)}(\bm y_{(C_i)})/ \lambda^{(S_i)}(\bm y_{S_i})\},
    \label{eq:lambda_dgm}
  \end{IEEEeqnarray}
  where
  $\lambda^{(A)}(\bm y_A) = \lVert \bm y_A\rVert^{-(|A|+1)} \, h^A(\bm
  y_A/\lVert \bm y_A\rVert)$, $\bm y_A \in \EE^A$ and
  $h^A(\bm \omega_A)$, $\bm \omega_A \in \SSS_{|A|-1}$ is the density
  of the spectral measure $H^A$ associated with the limit measure
  $\nu^A$. 
\end{proposition}
A proof of Proposition \ref{prop:maxstable} is given in Appendix
\ref{app:proof_prop_1}.\ Since the random subvector corresponding to any separator set must
have an identical distribution irrespective of the clique that
contains the separator, we can prove the following result on the
asymptotic properties of regularly varying random vectors
corresponding to these cliques.
\begin{proposition} 
  \label{prop:compatibility}
  Let $\bm X$ be a graphical model relative to a decomposable graph
  $\GG=(V,E)$.\ Then, for any two cliques $C$ and $C'$ for which
  $C\cap C'=S\neq 0$ such that the limit $\lim_{\level\rightarrow\infty}\level\, \PP\left(\bm X_{S}/\level
      \in [0,\bm x_S]^c\right)$ exists, then
  \begin{IEEEeqnarray}{rCl}
    &&\lim_{\level\rightarrow\infty}\level\, \PP\left(\bm X_{C}/\level
      \in B\right) = \lim_{\level\rightarrow\infty}\level\,
    \PP\left(\bm X_{C'}/\level \in B'\right),
    \label{eq:compatibility}
  \end{IEEEeqnarray}
  for any two sets 
  $B=\EE^{C \sm S}\times [0, \bm x_{S}]^c$ and
  $B'=\EE^{C' \sm S}\times [0, \bm x_{S}]^c$.
\end{proposition}
A proof of Proposition \ref{prop:compatibility} is given in Appendix
\ref{app:proof_prop_2}.  

Despite the mathematical convenience of the regular variation
assumption, it has the implication that the coefficient of extremal
dependence in Equation \ref{eq:chi} satisfies $\chi_{V}>0$ meaning
that the largest extremes of $\bm X$ occur simultaneously.\ Hence,
several existing approaches exclude random vectors that exhibit strong
dependence at finite levels, yet their largest extremes occur
independently in the limit, such as for Gaussian random vectors
arising from any Gaussian distribution with a strictly
positive-definite matrix.\ An approach that coalesces such distinct
types of extremal dependence in a unified treatment is missing.\ This
is the main focus of Section~\ref{sec:ce}.

\subsection{Conditional extreme value theory}
\label{sec:ce}
\color{black} For the remainder of the paper, we adopt an approach
based on conditional extreme value theory
\citep{hefftawn04,heffernan2007limit}.\ The key basis of this approach
rests on assuming that the marginals belong to the max-domain of
attraction of the Gumbel distribution, and are not regularly varying.\
For simplicity and clarity of exposition, we will make the following
assumption. 
\begin{description}[wide=0\parindent] \em
\item[Assumption $\AAA$.] The marginal distributions of $\bm X_E$ are
  all unit-exponential, that is,
  $\PP(X_{E,v} \leq x) = (1-\exp(-x))_+$, for any $v\in V$.
\end{description}
Conditional extreme value theory is a powerful approach that can be
used to characterize the distribution of a random vector conditionally
on the presence of an extreme event.\ In particular, conditional
extreme value theory presupposes that there exist norming functions
$\bm a^{(\condv)}_{V \sm \condv} : \RR \rightarrow \RR^{d-1}_+$,
$\bm b^{(\condv)}_{V \sm \condv}: \RR \rightarrow \RR^{d-1}_+$, and a
distribution function $G_{V\sm \condv}^{(v)}$ supported on
$\RR^{\lvert V\rvert-1}$ with non-degenerate margins such that, as
$\level\rightarrow \infty$,
\begin{equation}
  \lim_{\level\to\infty}\PP \left( X_{E,\condv}-\level > x_\condv,
    \frac{\bm X_{E, V \sm \condv} -\bm a^{(\condv)}_{V \sm \condv} (X_{E,\condv})}
    {\bm b^{(\condv)}_{V \sm \condv}(X_{E,\condv}) } \leq \bm {z}_{V \sm \condv} \ \Big | \ X_{E,\condv} >\level
  \right) = \exp(-x_{\condv})\, G^{(\condv)}_{V \sm \condv}\left ( \bm{z}_{V \sm \condv} \right),
  \label{eq:normdg}
\end{equation}
at continuity points
$(x_{\condv}, \bm z_{V \sm \condv}) \in \RR_+\times \RR^{|V|-1}$ of
the limit distribution function.\ 

To appreciate the generality of this approach, we remark that under
the prescribed marginal choice, when convergence \eqref{eq:normdg}
holds with $\bm a_{V\sm \condv}(\level)=\level\,\bm 1$ and
$\bm b_{V\sm \condv}(\level)=\bm 1$, then convergence
\eqref{eq:normdg} is equivalent to $\exp(\bm X_{E, \condv})$ being
multivariate regularly varying with a limit measure that places its
mass in the interior of $\mathbb{E}$ \citep{asenova2021extremes}.\
However, when convergence \eqref{eq:normdg} holds with
$\bm a_{V\sm \condv}^{(\condv)}(\level)/\level \to \bm \alpha_{V\sm
  \condv}^{(\condv)}\in [0,1)^{d-1}$ as $\level\to \infty$, then
affine normalization in exponential margins helps to reveal structure
which is hidden in subfaces of $\mathbb{E}$ under affine
normalization in regularly varying margins.\ Thus, conditional extreme
value theory has the additional key advantage that it can be adopted
to characterize extremes of both asymptotically dependent and
asymptotically independent random vectors.\


In the next section, we will analyze decomposable graphical models
through suitable conditions associated with the weak convergence of
the conditional distributions that define the Markov structure of the
joint distribution, specifically through the weak convergence of the
clique distribution conditionally on there being an extreme event in a
variable $X_{E,\condv}$ in this clique, and the weak convergence of
the transition probability
kernels.\ 
Thus, we will exploit the Markov property of decomposable graphical
models so that, under appropriate conditions, we can arrive at
convergence \eqref{eq:normdg}, that is, at affine normings for the
entire random vector $\bm X_{E, V \sm \condv}$ in terms of a variable
$X_{E,\condv}$.\ This approach of requiring specific limit behaviour
for a clique distribution and for the transition probability kernels
closely follows previous approaches on conditional extreme value
theory which focus on time-homogeneous Markov chains with continuous
state-space \citep{papastathopoulos_strokorb_tawn_butler_2017,
  papastathopoulos2023hidden}.\ Here, the emphasis is placed on
decomposable graphical models.
\color{black}

\section{Extremes of decomposable graphical models}
\label{sec:main}
\subsection{Decomposable tail graphical models}
\color{black} 
Without loss of generality, for any $\condv\in V$, we can choose a
clique $C_1$ in the running intersection property \eqref{eq:ri} such
that the vertex $\condv \in C_1$, see Remark \ref{rem:rip}.  In a
general decomposable graph, a reference vertex may belong to more than
one clique, as, for example, in the Goldner-Harary graph shown in
Figure \ref{fig:sepgh} where $\condv = 2$ is located in $6$ out of $8$
of the maximal cliques of the graph.  In this case, we choose an
arbitrary $C_1$ such that $v \in C_1$ and root the junction tree at
the clique $C_1$.\ As a result of the total order between cliques
induced by the running intersection property, we may informally
imagine that, after witnessing an extreme event $X_\condv > \level$,
where $\level$ is a large relatively to the marginal unit-exponential
distribution,
then the impact from the shock $X_{\condv}$ in the system, propagates
at any vertex $j\in V\sm \condv$ according to this total order.\ Our
next assumption guarantees that location and scale normings can be
found such that the conditional distribution of affinely renormalized
states $\bm X_{E, C_1 \sm \condv}$ given $X_{E,\condv} > \level$
converges weakly on $\RR^{| C_1 |-1}$ to a non-degenerate limit
distribution.
\begin{description}[wide=0\parindent] \em
\item [Assumption $\AA$.] There exist measurable functions
  $\bm a^{(\condv)}_{C_1 \sm \condv} : \RR \rightarrow \RR^{|C_1|-1}$
  and
  $\bm b^{(\condv)}_{C_1 \sm \condv }: \RR \rightarrow
  \RR^{|C_1|-1}_+$ , such that
  $\bm a^{(\condv)}_{C_1 \sm \condv }(\level) +\bm b^{(\condv)}_{C_1
    \sm \condv }(\level)\, \bm z_C \rightarrow \infty \, \bm
  1_{|C_1|-1} $ as $\level\rightarrow \infty$ for all
  $\bm z_{C_1 \sm \condv} \in \RR^{|C_1|-1}$, and a distribution
  function $G^{(\condv)}_{C_1 \sm \condv}$ supported on
  $ \RR^{|C_1|-1}$ with non-degenerate marginal distributions, such
  that as $\level\rightarrow \infty$,
  \begin{equation}
    \PP\left(\frac{\bm X_{E,C_1 \sm \condv } -\bm a^{(\condv)}_{C_1 \sm \condv } (X_{E,\condv})}
      {\bm b^{(\condv)}_{C_1 \sm \condv}(X_{E,\condv}) } \leq
      \bm {z}_{C_1 \sm \condv}~\Big |~ X_{E,\condv} = \level \right)
    \wk G^{(\condv)}_{C_1  \sm \condv }\left ( \bm{z}_{C_1 \sm \condv } \right),
    \qquad \bm{z}_{C_1 \sm \condv} \in \RR^{|C_1|-1}.
    \label{eqn:a1}
  \end{equation}
\end{description}
Now, as we transit away from the initial clique, we need to ensure the
weak convergence of the conditional distribution of variables indexed
by cliques given variables indexed by
separators. 
This is the role of our next assumption which requires 
that conditionally on $\bm X_{S_{\ttime}}$ growing at an appropriate
speed, then location and scale functionals can be found so that the
conditional distribution of the $\bm X_{C_\ttime\sm S_\ttime}$
affinely normalized in terms of $X_{S_\ttime}$ converges weakly to a
non-degenerate distribution.\ This approach is similar to that in
\cite{papastathopoulos2023hidden} for the special case of a
higher-order Markov chain, but here it is extended to the case of
decomposable graphical models.  \color{black}
\begin{description}[wide=0\parindent]
\item[Assumption $\AB$.] \em
\begin{enumerate}[wide=0\parindent] 
\item[$(i)$] There exist vector-valued measurable functions
  $\bm a^{(\condv)}_{V \sm \condv} : \RR_+ \rightarrow \RR^{|V|-1}_+$,
  and
  $\bm b^{(\condv)}_{V \sm \condv}: \RR_+ \rightarrow \RR^{|V|-1}_+$
  such that
  $\bm a^{(\condv)}_{V \sm \condv}(\level)+ \bm b^{(\condv)}_{V \sm
    \condv}(\level) \bm z \rightarrow \infty \bm 1_{|V|-1}$ as
  $\level \rightarrow \infty$ for all $\bm z \in \RR^{|V|-1}$.
  Furthermore, $\bm a^{(\condv)}_{C_1 \sm \condv}$,
  $\bm b^{(\condv)}_{C_1 \sm \condv}$ are the norming functions
  specified in Assumption $\AA$.  For $\ttime = 2, \dots, N$, there exist
  vector-valued measurable functions
  $\bm a^{(S_{\ttime})}_{C_{\ttime} \sm S_{\ttime}}:
  \RR^{|S_{\ttime}|}_+ \rightarrow \RR^{|C_{\ttime} \sm
    S_{\ttime}|}_+$ and 
  $\bm b^{(S_{\ttime})}_{C_{\ttime} \sm S_{\ttime}}:
  \RR^{|S_{\ttime}|}_+ \rightarrow \RR^{|C_{\ttime} \sm
    S_{\ttime}|}_+$, and continuous update functions
  $\bm \psi^{(S_{\ttime})}_{C_{\ttime} \sm S_{\ttime}} :
  \RR^{|S_{\ttime}|} \rightarrow \RR^{|C_{\ttime} \sm {S_{\ttime}|}} $
  and
  $\bm \phi^{(S_{\ttime})}_{C_{\ttime} \sm S_{\ttime}} :
  \RR^{|S_{\ttime}|} \rightarrow \RR_+^{{|C_{\ttime}} \sm S_{\ttime}|}
  $, such that for all $\bm z_{S_{\ttime}} \in \RR^{|S_{\ttime}|}$,
    
    \begin{equation}
      \bm \psi^{(S_{\ttime})}_{C_{\ttime} \sm S_{\ttime}} (\bm z_{S_{\ttime}}) = \lim_{\level \to
        \infty} \frac{\bm a^{(S_{\ttime})}_{C_{\ttime} \sm S_{\ttime}} \{\bm T_{S_{\ttime}} (\bm z_{S_{\ttime}},
        \level) \} -
        \bm a^{(\condv)}_{C_{\ttime} \sm S_{\ttime}}(\level)}{\bm b_{C_{\ttime} \sm S_{\ttime}}^{(\condv)} (\level)} \quad
      \text{and} \quad \bm \phi^{(S_{\ttime})}_{C_{\ttime} \sm S_{\ttime}} (\bm z_{S_{\ttime}}) = \lim_{\level \to
        \infty} \frac{ \bm b^{(S_{\ttime})}_{C_{\ttime} \sm S_{\ttime}} \{ \bm T_{S_{\ttime}} (\bm z_{S_{\ttime}},
        \level) \} }{\bm b_{C_{\ttime} \sm S_{\ttime}}^{(\condv)} (\level)} ,
    \label{eq:up2}
  \end{equation}
   where $\bm T^(\condv)_{S_{\ttime}} (\bm z_{S_{\ttime}}, \level) = \bm a^{(\condv)}_{S_{\ttime}}(\level) + \bm
  b^{(\condv)}_{S_{\ttime}}(\level) \bm z_{S_{\ttime}}$ and such that the remainder terms,
  \begin{IEEEeqnarray}{rCl}
    \bm  A^{(S_{\ttime})}_{C_{\ttime} \sm S_{\ttime}}(\bm z_{S_{\ttime}}, \level) &=& \frac{\bm a^{(\condv)}_{C_{\ttime} \sm S_{\ttime}}(\level) - \bm a^{(S_{\ttime})}_{C_{\ttime} \sm S_{\ttime}}\{ T^{(\condv)}_{S_{\ttime}} (\bm z_{S_{\ttime}},         \level) \}  + \bm b_{C_{\ttime} \sm S_{\ttime}}^{(\condv)} (\level) \bm  \psi^{(S_{\ttime})}_{C_{\ttime} \sm S_{\ttime}} (\bm z_{S_{\ttime}})}{\bm b^{(S_{\ttime})}_{C_{\ttime} \sm S_{\ttime}} \{ T^{(\condv)}_{S_{\ttime}} (\bm z_{S_{\ttime}},         \level) \} } \ \text{and} \nonumber \\
    \bm B^{(S_{\ttime})}_{C_{\ttime} \sm S_{\ttime}}(\bm
    z_{S_{\ttime}}, \level) &=& 1- \frac{\bm b^{(\condv)}_{C_{\ttime}
        \sm S_{\ttime}} (\level) \bm \phi^{(S_{\ttime})}_{C_{\ttime}
        \sm S_{\ttime}} (\bm z_{S_{\ttime}})}{\bm
      b^{(S_{\ttime})}_{C_{\ttime} \sm S_{\ttime}}\{ \bm T^{(\condv)}_{S_{\ttime}}
      (\bm z_{S_{\ttime}}, \level) \} },
      \label{eq:remainder}
  \end{IEEEeqnarray}
  converge to zero as $\level \to \infty $, uniformly in compact sets
  of $\bm z_{S_{\ttime}} \in \RR^{|S_{\ttime}|}$.
\item[(ii)] For $\ttime = 2, \dots, N$, there exist distributions
  $G_{C_{\ttime} \sm S_{\ttime}}^{(S_{\ttime})}$ supported on
  $\RR^{|C_{\ttime} \sm S_{\ttime}|}$, with non-degenerate marginal
  distributions such that as $\level \rightarrow \infty$,
  \begin{equation} \PP \left (\frac{ \bm X_{E,C_{\ttime} \sm
          S_{\ttime}}- \bm a^{(S_{\ttime})}_{C_{\ttime} \sm
          S_{\ttime}}(\bm X_{E,S_{\ttime}})}{\bm
        b^{(S_{\ttime})}_{C_{\ttime} \sm S_{\ttime}}(\bm
        X_{E,S_{\ttime}})} \leq \bm y_{C_{\ttime} \sm S_{\ttime}}
      ~\Big |~ \frac{\bm X_{E,S_{\ttime}} - \bm
        a^{(\condv)}_{S_{\ttime}}(\level)} {\bm
        b^{(\condv)}_{S_{\ttime}}(\level)} = \bm z_{S_{\ttime}}
    \right) \wk G_{C_{\ttime} \sm S_{\ttime}}^{(S_{\ttime})}\left (
      \bm y_{C_{\ttime} \sm S_{\ttime}}\right),
      \label{eqn:a2}
    \end{equation}
    uniformly in compact sets in the
    variable $\bm z_{S_{\ttime}} \in \RR^{|S_\ttime|}$.    
  \end{enumerate}
\end{description}
 
\begin{remark} \normalfont It is convenient for the remainder of the
  paper to define $S_1 = \{v \}$, and to set
  $\bm \psi^{(S_{1})}_{C_{1} \sm S_{1}} (\bm x_{S_1}) =0$ and
  $\bm \varphi^{(S_{1})}_{C_{1} \sm S_{1}}(\bm x_{S_1}) =1$, for all
  $\bm x_{S_1} \in \RR^{|S_1|}$.
  \label{rem:s1}
\end{remark}
With these assumptions in place, we are in position to state our first
main
theorem.
\begin{theorem} \label{thm:main_theorem} Let $\bm X_E$ be graphical
  model relative to a connected, undirected, decomposable graph $\GG$,
  \ that satisfies Assumptions $\AAA$, $\AA$ and $\AB$.\ Then, as
  $\level\rightarrow \infty$, 
  \begin{equation}
    \left({X_{E,\condv}-\level} , \frac{\bm X_{E,V \sm \condv}-\bm a^{(\condv)}_{V \sm \condv}(X_{E,\condv})}{\bm b^{(\condv)}_{V \sm \condv}(X_{E,\condv})}\right) ~\Bigg |~ \{X_{E,\condv} >\level\}  \cind ( E_{\condv},  \bm Z^{(\condv)}_{V \sm \condv}  ),
    \label{eq:th11}
  \end{equation}
  where
  \begin{enumerate}[wide=0\parindent] 
  \item[$(i)$] $E_{\condv}$ is a unit exponential random variable that
    is independent of $\bm Z^{(\condv)}_{V \sm \condv}$,
     
  \item[$(ii)$] for $\ttime = 1, \dots, N$,
  \begin{equation}
    \bm Z^{(\condv)}_{C_{\ttime}  \sm S_{\ttime}} = \bm \psi^{(S_{\ttime})}_{C_{\ttime}  \sm S_{\ttime}} (\bm Z^{(\condv)}_{S_{\ttime}}) +  \bm \phi^{(S_{\ttime})}_{C_{\ttime}  \sm S_{\ttime}} (\bm Z^{(\condv)}_{S_{\ttime}}) \,\bm \varepsilon^{(\condv)}_{C_{\ttime}  \sm S_{\ttime}},
    \label{eq:tailgraph}
  \end{equation}
  where $Z^{(\condv)}_{\condv}=0$ a.s.,
  $\bm \varepsilon^{(\condv)}_{C_{\ttime} \sm S_{\ttime}} \sim
  G^{(S_{\ttime})}_{C_{\ttime} \sm S_{\ttime}}$ for
  $\ttime = 1, \dots, N$, and the random vectors
  $\bm \varepsilon^{(\condv)}_{C_{\ttime_1} \sm S_{\ttime_1}}$ and
  $\bm \varepsilon^{(\condv)}_{C_{\ttime_2} \sm S_{\ttime_2}}$ are
  independent when $\ttime_1 \neq \ttime_2$.
  \end{enumerate}
 \label{theorem:main}
\end{theorem}

\subsection{Tail noise for block graphs}
\color{black} Here, we restrict consideration to graphical models with
conditional independence relationships described by a block graph,
that is, by a chordal graph $\GG$ where all separating sets $S_\ttime$
in the running intersection \eqref{eq:ri} satisfy $|S_\ttime|=1$.\ We
focus mainly on block graphical models that are specified via their
clique distributions $F_{\bm X_{E,C_\ttime}}$, where the cliques
exhibit different extremal dependence structures.\ In particular, we
are interested in graphical models where some cliques exhibit
asymptotic dependence whilst others exhibit asymptotic independence.\
Although such graphical models are straightforward to construct, they
have received no attention in the literature of graphical modelling of
extremes.\

To understand the behaviour of extreme events of such block graphical
models, we first note that Assumption $\AB$ may be violated when the
norming functions in one clique grow on a different scale to those in
its neighbouring clique.\ 
We demonstrate that in cases where Assumption $\AB$ breaks, we can
still prove convergence with norming functions of the variables
associated with the separators, rather than with norming functions of
the variable that is conditioned to grow large.\ This change of
norming functions permits us This allows a wide range of combinations
of different clique distributions to be considered.
\color{black}
\begin{description}[wide=0\parindent] \em
\item [Assumption $\BA$.] Let $\GG$ be a connected block graph.\ For all $C_{\ttime}\in \mathcal{M}(\GG)$ and for any vertex $v_{\ttime}\in
  C_{\ttime}$,\ there exist measurable functions
  $\bm a^{(\condv_{\ttime})}_{C_{\ttime} \sm v_{\ttime}} : \RR \rightarrow
  \RR^{|C_{\ttime}|-1}$ and
  $\bm b^{(\condv_{\ttime})}_{C_{\ttime} \sm v_{\ttime} }: \RR \rightarrow
  \RR^{|C_{\ttime}|-1}_+$,\ such that
  $\bm a^{(\condv_{\ttime})}_{C_{\ttime} \sm v_{\ttime} }(\level) +\bm
  b^{(\condv_{\ttime})}_{C_{\ttime} \sm v_{\ttime} }(\level)\, \bm y_C
  \rightarrow \infty \, \bm 1_{|C_{\ttime}|-1} $ as
  $\level\rightarrow \infty$ for all
  $\bm y_{C_{\ttime} \sm v_{\ttime}} \in \RR^{|C_{\ttime}|-1}$, and a
  distribution function $G^{(\condv_{\ttime})}_{C_{\ttime} \sm \condv_{\ttime}}$
  supported on $ \RR^{|C_{\ttime}|-1}$ with non-degenerate marginal
  distributions, such that as $\level\rightarrow \infty$,
  \begin{equation}
    \PP\left(\frac{\bm X_{E,C_{\ttime} \sm v_{\ttime} } -\bm a^{(v_{\ttime})}_{C_{\ttime} \sm v_{\ttime} } (X_{E,v_{\ttime}})}
      {\bm b^{(v_{\ttime})}_{C_{\ttime} \sm v_{\ttime}}(X_{E,v_{\ttime}}) } \leq
      \bm {z}_{C_{\ttime} \sm v_{\ttime}}~\bigg |~ \{X_{E,\condv_{\ttime}}\} = \level \right)
    \cind G^{(v_{\ttime})}_{C_{\ttime}  \sm v_{\ttime} }\left ( \bm{z}_{C_{\ttime} \sm v_{\ttime} } \right),
    \quad \bm{z}_{C_{\ttime} \sm v_{\ttime}} \in \RR^{|C_{\ttime}|-1}.
    \label{eqn:b1}
  \end{equation}
\end{description}
With these assumptions we apply a normaliztion in terms of
$\bm X_{S_{\ttime}}$ and obtain a weak limit of the normalised random
vector as $t \rightarrow \infty$.  We use the term tail noise to
denote this weak limit.
\begin{theorem} 
  \label{prop:block}
  Let $\bm X_{E}$ be a graphical model relative to a connected
  undirected block graph $\GG$.\ If $\bm X_E$ satisfies Assumptions
  $\AAA$ and $\BA$, then as $\level \to \infty$,
\begin{equation}
  \left ( \frac{ \bm X_{E,C_{\ttime} \sm S_\ttime } -
      \bm a^{(S_\ttime)}_{C_{\ttime} \sm S_\ttime } (X_{E,S_\ttime}) }
    {\bm b^{(S_\ttime)}_{C_{\ttime} \sm S_\ttime}(X_{E,S_\ttime})}\,:\,
    \ttime =1,\dots,N \right) 
  \ \bigg | \ \{X_{E,S_1} > \thres\}  \cind
  \left (\bm Z_{C_\ttime \sm S_\ttime}^{(S_\ttime)}\,:\,\ttime=1,\dots,N \right),  
     \label{eq:th2}
\end{equation}
where
$\bm Z_{C_{\ttime} \sm S_{\ttime}}^{(S_\ttime)} \sim G_{C_{\ttime} \sm
  S_{\ttime}}^{(S_{\ttime})}(\bm y_{C_{\ttime} \sm S_{\ttime}})$ and
$\bm Z_{C_{\ttime} \sm S_{\ttime}}^{(S_\ttime)} \bigCI \bm Z_{C_{j}
  \sm S_{j}}^{(S_j)}$ for $i \neq j$.
\end{theorem}
\color{black} Here, we remark that the choice of the term tail noise
is made to reflect the independence of
$\bm Z_{C_{\ttime} \sm S_{\ttime}}^{(S_{\ttime})}$ and
$\bm Z_{C_{j} \sm S_{j}}^{(S_{j})}$ for all $\ttime \neq j$. The
elements of each $\bm Z_{C_{\ttime}\sm S_{\ttime}}^{S_{\ttime}} $,
$i=1,\dots,N$, need not be
independent.
\section{Examples}
\label{sec:examples}

\subsection{Asymptotic behaviour of transition probability kernels of
  multivariate max-stable distributions}
We assume that $\Lambda^{(C)}$, the exponent measure in
\eqref{eq:expfn} obtained from the distribution of the random variable
$\bm X_{E,C}$, has mixed derivatives of all orders, this implies that
the joint distribution of $\bm X_C$ admits a joint probability density
function \citep{coles1991modelling}.\ 
Let $\Pi_J$ be the set of all partitions of
$J\subseteq C$ and
$\Lambda_J^{(C)}(\bm y) = \partial^{|J|} \Lambda^{(C)}(\bm y)/\prod_{j \in J}
\partial x_j$.\ Let
$\Lambda^{(S)}(\bm y_S) := \lim_{\bm y_{C \sm S} \rightarrow \infty\bm
  1_{|C \sm S|}} \Lambda^{(C)}( \bm y)$ and
$\Lambda_J^{(S)}( \bm y_{S}) = \partial^{|J|}\Lambda^{(S)}(\bm y_{S})
/\prod_{j \in J} \partial x_j $ for any non-empty $J \subseteq S$ and.\ The transition probability kernel is equal
to 
\begin{equation}
  \PP\left (\bm X_{E, C \sm S }\leq \bm x_{C \sm S } \mid \bm X_{E,S}= \bm x_S \right ) = \frac{\sum_{p \in \Pi_S} (-1)^{|p|} \prod_{J \in p} \Lambda_J^{(C)}(\bm y)}{ \sum_{p \in \Pi_S} (-1)^{|p|} \prod_{J \in p} \Lambda_J^{(S)}(\bm y_S)} \exp \left \{ \Lambda^{(S)}(\bm y_S) - \Lambda^{(C)}(\bm y) \right \},
  \label{eq:MS1}
\end{equation}
where $\bm y = -1/\log\{1-\exp(- \bm x)\}$, $\bm x \in \RR_+^d$.

We seek a functional $\bm a^{(S)}_{C \sm S}(\bm x_S)$ such that the
probability
$\PP (\bm X_{E,C \sm S }\leq \bm a^{(S)}_{C \sm S }(\bm X_{E,S}) \mid
\bm X_{E,S}= \level\,\bm 1_{|S|} + \bm z_S )$ converges to a value in
the interval $(0,1)$.\ See \cite{papastathopoulos2023hidden} for a
similar analysis.\ Taking equation \eqref{eq:MS1} and setting
$\bm x_S= \level \bm 1_{|S|} + \bm z_S$ and
$\bm x_{C \sm S } = a^{(S)}_{C \sm S }(\level\,\bm 1_{|S|} + \bm
z_S)$, we have that as $\level \rightarrow \infty$, equation
\eqref{eq:MS1} reduces to,
  \begin{equation*}
    \PP\left (\bm X_{E,C \sm S }\leq \bm a^{(S)}_{C \sm S }(\bm X_{E,S}) \mid \bm X_{E,S}= \level\,\bm 1_{|S|} + \bm z_S \right ) = \{ \Lambda^{(C)}_S(\bm y)/\Lambda^{(S)}_S(\bm y_S) \}
    (1+o(1)) \quad \text{as $\level\to\infty$},
  \end{equation*}
  and for any continuous functional $\bm a^{(S)}_{C \sm S }$ that
  satisfies
  \begin{equation}
    \exp [\bm a^{(S)}_{C \sm S }\{\log (\level \, \bm x_S)\}] = \level\,\exp [\bm
    a^{(S)}_{C \sm S }\{\log (\bm x_S)\}],
    \label{eq:ms_a_hom}
  \end{equation}
  we get
  \begin{equation}
    \lim_{\level \rightarrow \infty} \PP\left (\bm X_{E, C \sm S }\leq \bm a^{(S)}_{C \sm S }(\bm X_S) \mid \bm X_{E,S}= \level\,\bm 1_S + \bm z_S \right ) = \frac{\Lambda^{(C)}_S \{ \exp( \bm z_S), \exp( \bm a^{(S)}_{C \sm S }(\bm z_S)\}}{\Lambda^{(S)}_S \{ \exp(\bm z_S)) \}}.
    \label{eq:maxstableres}
  \end{equation}
\color{black}\subsection{Decomposable graphical models with
  H\"{u}sler--Reiss max-stable cliques}
\label{eq:HR_cliques_example} Let
$\bm X_E = F_E^{-1}\{ F_F(\bm X_F)\}$ where $\bm X_F$ is a graphical
model with unit Fr\'{e}chet marginal distributions.\ Here we consider
graphical models with clique distributions being max-stable, that is,
$F_{\bm X_{F,C}}(\bm x_C) =\exp\{-\Lambda^{(C)}(\bm x_C) \}$ for
$\bm X_C \in \RR_+^{|C|}$, where $\Lambda^{(C)}(\bm x_C)$ denotes the
exponent measure of the H\"{u}sler--Reiss max-stable distribution,
which is given by \eqref{eq:HR_cdist} .\ The H\"{u}sler--Reiss
max-stable distribution is parameterized by a variogram matrix
$\bm \Gamma_C =(\Gamma_{ij})_{i,j\in C}\in \RR^{|C|\times|C|}$
\citep{kabluchko2009stationary,engelke2020graphical}, that is, by a
strictly conditionally negative definite matrix satisfying
$\text{diag}(\bm \Gamma_C)=0$ and
$\bm a^\top \,\bm \Gamma_C \,\bm a < 0$ for any
$\bm a=(a_i)_{i\in C}\in\RR^{|C|}$ with $\sum_{i \in C}
a_i=0$.\ 
The exponent measure $\Lambda^{(C)}$ is defined by 

\begin{equation}
  \Lambda^{(C)}(\bm x_C; \bm \Gamma_C) = \sum_{c \in C} \frac{1}{x_{c}}
  \Phi_{|C|-1} \left \{ \log\left(\frac{\bm x_{C\sm c}}{x_{c}}\right) + \frac{\bm \Gamma_{c, C\sm c}}{2};  \bm 0_{C},
    \bm \Sigma_{C\sm c}^{(c)} \right \},
  \label{eq:HR_cdist}
\end{equation}
where
$2\, \bm \Sigma_{C\sm c}^{(c)}= \bm \Gamma_{C\sm c,
  c} \, \bm 1^\top + (\bm \Gamma_{C\sm c, c}\, \bm
1^\top)^\top + \bm \Gamma_{C\sm c}$ and
$\Phi_{|C|-1}(\cdot\, ;\, \bm 0_{C}, \bm \Sigma_{C\sm
  c}^{(c)})$ denotes the cumulative distribution function of
the multivariate normal with mean vector $\bm 0_C$ and covariance
matrix $\bm \Sigma_{C\sm c}^{(c)}$.\ Assumption $\AA$ is
satisfied since, from \cite{engeetal15} (equation 47),  for all $v \in V$ 
and for any clique $C_1$ such that $v\in C_1$, then 
\begin{equation*}
  G^{(\condv)}_{C_1 \sm \condv}(\bm z_{C_1 \sm \condv}) = \Phi_{|C_1|-1}\left  (\bm z_{C_1 \sm \condv}, -\text{diag}(\bm \Sigma_{C_1}^{(\condv)})/2,\ \bm \Sigma_{C_1}^{(\condv)} \right ).
\end{equation*}
The limit distribution for the renormalized transition probability kernel described in Assumption $\AB$ is
stated in Lemma \ref{lemma:hra2}.
\begin{lemma}
  For $i=1,\dots, N$, let $s_1=v $ and $s_{\ttime}=v$ if $v \in S_{\ttime}$.\ Otherwise if $\condv \notin S_{\ttime}$ let
  $s_{\ttime}$ be an arbitrary vertex in $S_{\ttime}$. Let
  $\bm Q^{(S_{\ttime})}_{C_{\ttime}} = \big( \bm
  \Sigma^{(S_{\ttime})}_{C_{\ttime} \setminus s_{\ttime}} \big )^{-1}$
  and define the matrix
  $\widetilde{\bm Q}_{C_{\ttime}\sm S_{\ttime}, S_{\ttime}}$ by
  $\widetilde{\bm Q}_{C\sm S,S\sm s}=\bm Q^{(s_{\ttime})}_{C\sm S,S \sm
    s}$ and \
  $\widetilde{\bm Q}_{C_{\ttime}\sm S_{\ttime}, s_{\ttime}}= -\bm
  Q^{(s_{\ttime})}_{C_{\ttime}\sm S_{\ttime}, C_{\ttime}\sm
    s_{\ttime}}\,\bm 1_{|C_{\ttime}\sm s_{\ttime}|}$.\ With the choice of
  $\bm a^{(S_{\ttime})}_{C_{\ttime} \sm S_{\ttime} }(\bm x_{S_{\ttime}})
  = - \big ( \bm Q^{(s_{\ttime})}_{C_{\ttime} \sm S_{\ttime}} \big
  )^{-1} \widetilde{\bm Q}_{C_{\ttime}\sm S_{\ttime},S_{\ttime}}\, \bm
  x_{S_{\ttime}}$ and of
  $\bm b^{(S_{\ttime})}_{C_{\ttime} \sm S_{\ttime}}(\bm x_{S_{\ttime}})=
  1$, then Assumption $\AB$ is satisfied with
  \begin{equation}
    G^{(S_{\ttime})}_{C_{\ttime}\sm S_{\ttime}} (\bm y_{C_{\ttime}\sm
      S_{\ttime}})=\Phi_{|C_{\ttime}\sm S_{\ttime}|} \left ( \bm
      y_{C_{\ttime} \sm S_{\ttime}}\,;\,\bm m_{C_{\ttime} \setminus
        S_{\ttime}}, (\bm Q^{(s_{\ttime})}_{C_{\ttime} \sm
        S_{\ttime}})^{-1}\right), \quad \bm
    \psi^{(S_{\ttime})}_{C_{\ttime} \sm S_{\ttime}} (\bm
    x_{S_{\ttime}})= \bm a_{C_{\ttime}\sm S_{\ttime}}^{(S_{\ttime})}(\bm
    x_{S_{\ttime}}), \quad\bm \varphi^{(S_{\ttime})}_{C_{\ttime} \sm
      S_{\ttime}} (\bm x_{S_{\ttime}})=1,
    \label{eq:l1}
  \end{equation}
  where
  $ \bm m_{C_{\ttime} \setminus S_{\ttime}} = - \big(\bm
  Q^{(s_{\ttime})}_{ C_{\ttime}\sm S_{\ttime}}\big)^{-1} \bm
  Q^{(s_{\ttime})}_{ C_{\ttime}\sm S_{\ttime}, C_{\ttime}\sm s_{\ttime}}
  \, \bm \Gamma_{C_{\ttime}\sm s_{\ttime}}/2$.
  \label{lemma:hra2}
\end{lemma}   
A proof of Lemma \ref{lemma:hra2} is included in Appendix C.1.\ The
expressions
$\bm a_{C_{\ttime}\sm S_{\ttime}}^{(S_{\ttime})}, \widetilde{\bm
  Q}_{C_{{\ttime}}\sm S_{{\ttime}}, S_{{\ttime}}}$ and
$\bm m_{C_{\ttime} \setminus S_{\ttime}}$ in Lemma \ref{lemma:hra2} do
not depend on the choice of $s_{\ttime} \in S_{{\ttime}}$ due to Lemma
1 of \cite{engelke2020graphical}, which asserts that
$\bm Q_{C_{\ttime}\sm S_{\ttime}}^{(s_{\ttime})}$ is independent of
$s_{\ttime}$ when $s_i\in S_i$.\ Hence, from Theorem
\ref{theorem:main}, $\bm Z_{S_1}^{(S_1)} =0$ a.s.,
$\bm Z_{C_1\sm S_1}^{(S_1)}\sim G_{C_1\sm S_1}^{(S_1)}$ and
\begin{IEEEeqnarray}{lClrCr}
  \bm Z^{(S_1)}_{C_{\ttime} \sm S_{\ttime}} &=& -(\bm
  Q^{(s_{\ttime})}_{C_{\ttime}\sm S_{\ttime}})^{-1} \widetilde{\bm
    Q}_{C_{\ttime}\sm S_{\ttime}, S_{\ttime}} \, \bm
  Z^{(S_1)}_{S_{\ttime}} + \bm \varepsilon_{C_{\ttime} \sm
    S_{\ttime}}, \quad \ttime =1, \dots N,
  \label{eqn:tgexp}
\end{IEEEeqnarray}
where $\bm \varepsilon_{C_{\ttime} \sm S_{\ttime}}\sim
G^{(S_{\ttime})}_{C_{\ttime}\sm S_{\ttime}}$ with $G^{(S_{\ttime})}_{C_{\ttime}\sm S_{\ttime}}$ given by expression \eqref{eq:l1}.
To fully describe the distribution of $\bm Z_{V\sm \condv}^{(\condv)}$
we introduce some additional notation.\ For any matrix
$\bm M = [M_{n m}]_{n\in V_1, m\in V_2 }$, indexed by the vertex
subsets $V_1 \subset V\sm \condv$ and $V_2 \subset V\sm \condv$ ,
define the $(|V|-1) \times (|V|-1) $-dimensional matrix
$\breve{\bm M} = [M_{n m}]_{n, m\in V \bm \condv }$ such that $\breve{M}_{n, m} = {M}_{n, m}$ for $n \in V_1$ and $m \in V_2$ and $\breve{M}_{n, m} = 0$ otherwise.

Applying results concerning the inverse of sparse precision matrices
\citep{JOHNSONlocalinversion}, we derive expressions for the mean and
precision of the Gaussian random vector
$\bm Z^{(\condv)}_{V \sm
  \condv}$ , the weak limit of the random vector
 \begin{equation}
 \left (\bm X_{E, V \sm \condv} - t \bm 1_{V \sm \condv} \ | \ X_{E,v} =t  \right ) \quad \text{as} \quad t \rightarrow \infty,
 \label{eq:stateHR}
\end{equation} 
\begin{proposition}
  \label{prop:HR} \normalfont Let $\bm Z^{(\condv)}_{V\sm \condv}$ be
  the weak limit in \eqref{eq:stateHR} then
  $\bm Z^{(\condv)}_{V\sm \condv}$ has distribution
  \smash{$\mathcal{N}_{d-1}\big ( \bm \mu^{v}_{V\sm v}, \bm
    \Sigma_{V\sm v}^{v} \big )$}, where the mean vector
  \smash{$\bm \mu^{v}_{V\sm v}$} is defined iteratively via 
  \begin{IEEEeqnarray*}{rCl}
    \bm \mu_{C_{\ttime}}^{(\condv)} &=& -\text{diag}\big (\bm
    \Sigma_{C_{\ttime}}^{(\condv)} \big )/2,
    \quad \text{when $C_{\ttime}\cap\{v\}\neq \emptyset$, and}  \\
    \bm \mu^{(\condv)}_{C_{\ttime} \sm S_{\ttime}} &=& -(\bm
    Q^{(s_{\ttime})}_{C_{\ttime}\sm S_{\ttime}})^{-1} \widetilde{\bm
      Q}_{C_{\ttime}\sm S_{\ttime}, S_{\ttime}} \, \bm
    \mu^{(\condv)}_{S_{\ttime}} + \bm m_{C_{\ttime} \setminus S_{\ttime}}, \quad
    \text{otherwise},
  \end{IEEEeqnarray*}
  and the inverse of the covariance matrix
  $\big ( \bm Q_{V \sm \condv}^{v} \big ):= \big(\bm \Sigma_{V\sm
    v}^{v}\big)^{-1}$ satisfies
  \[
    \bm Q_{V \sm \condv}^{v} = \sum_{\ttime=1}^N \breve { \bm
      Q}^{\ttime} - \sum_{\ttime=2}^N \breve{\bm P}^{\ttime-1} ,
\]
  where ${ \bm Q}^{\ttime}$ and ${ \bm P}^{\ttime-1}$ are defined iteratively through,
  \begin{IEEEeqnarray*}{lCl}
    \bm Q^{\ttime} &=& \left [
      \begin{array}{cc}
        \bm Q^{(s_{\ttime})}_{C_{\ttime} \sm S_{\ttime}} +   \bm A_{C_{\ttime} \sm S_{\ttime}, S_{\ttime}} {\bm P}^{\ttime-1} \Big (\bm A_{C_{\ttime}\sm S_{\ttime}, S_{\ttime}} \Big )^T & -  \bm A_{C_{\ttime}\sm S_{\ttime}, S_{\ttime}} {\bm P}^{\ttime-1}\\[8pt]
        -{\bm P}^{\ttime-1}\Big ( \bm A_{C_{\ttime}\sm S_{\ttime}, S_{\ttime}} \Big  ) ^T &  {\bm P}^{\ttime-1}
      \end{array} \right ]  \quad \text{when $C_{\ttime} \cap \{v\}= \emptyset$} , 
       \end{IEEEeqnarray*}
       where
  \begin{equation}      
      \bm A_{C_{\ttime}\sm S_{\ttime}, S_{\ttime}} = -\Big (\bm
  Q^{(s_{\ttime})}_{C_{\ttime}\sm S_{\ttime}}\Big )^{-1} \widetilde{\bm
    Q}_{C_{\ttime}\sm S_{\ttime},S_{\ttime}}
      \quad \text{and} \quad
    {\bm P}^{\ttime-1} = \left [ \left ( \left \{ \bm Q^{\ttime-1} \right \} ^{-1} \right ) _{S_{\ttime},S_{\ttime} } \right ]^{-1}.
      \end{equation}
\end{proposition}
\subsection{Decomposable meta-Gaussian graphical models}
\label{sec:gaussian}
 \color{black} Let $\bm X_E=F_{E}^{-1}(\Phi(\bm X_N))$ where
  $\bm X_N\sim \mathcal{N}_d(\bm 0, \bm R)$ is a zero-mean Gaussian
  graphical model on a decomposable graph $\GG$ with positive-definite
  covariance matrix $\bm R = (\rho_{ij})_{ij=1}^{d}$, $\rho_{ij}>0$,
  $\rho_{ii}=1$.\ Let
  $\bm \rho_{A,\condv}= (\rho_{i\condv}\,:\,i\in A)$, for any
  $A\subseteq V$, and write $\bm Q=(q_{ij})_{ij=1}^d=\bm R^{-1}$ for
  the precision matrix of the multivariate normal distribution.\
  \color{black} \cite{hefftawn04} showed that Assumption $\AA$ holds
  true with
  $\bm a^{(\condv)}_{C}(\level) = \bm \rho_{C \sm \condv}^2\,\level$ and
  $\bm b^{(\condv)}_{C}(\level) = \level^{1/2} \bm 1_{|C|}$.\ The limit
  distribution in expression \eqref{eqn:a1} is
  $G_{C}^{(\condv)}(\bm x_C)=\Phi_{|C|}(\bm x_C; \bm 0, \bm
  R_{\condv}^{C})$, $\bm x\in\RR^{|C}$, with
  $\bm R_{C}^{(\condv)} = \left( 2 \rho_{i\condv}
    \rho_{j\condv}(\rho_{ij} - \rho_{i\condv} \rho_{j\condv} )
  \right)_{i,j \in C }$.\ 
  The results from Example 1 of \cite{papastathopoulos2023hidden} can
  be generalised to show that, for the choice of functionals
  \begin{equation*} \bm a^{(S)}_{C\sm S}(\bm x_S) = \left(-\bm
      Q_{C\sm S}^{-1}\bm Q_{C\sm S, S} \,
      |\bm x_S|^{1/2}\right)^2\quad\text{and}\quad \bm b^{(S)}_{C\sm
      S}(\bm x_S) = \bm a^{(S)}_{C\sm S}(\bm x_S)^{1/2}, \qquad \bm x_S
    \in \RR^S,
  \end{equation*}
  Assumption $\AB$ is satisfied with limit distribution in expression
  \eqref{eqn:a2} equal to
  $G_{C\sm S}^{(S)}(\bm x_S) = \Phi(\bm x_S ;\bm 0, 2\, \bm Q_{C\sm
    S}^{-1}\color{black})$, and
  $\bm \psi_{C_{} \sm S_{}}^{S_{}} (\bm x)=\bm J_{\bm a_{C \sm
      S}^{(S)}}(\bm \rho_{S,\condv}^2)\, \bm x_S
  $ and
  $\bm \phi^{(S_{\ttime})}_{C_{\ttime} \sm S_{\ttime}} (\bm x_S)=\bm
  J_{\bm a_{C \sm S}^{(S)}}(\bm
  \rho_{S,\condv}^2)^{1/2}
  $, where
  $\bm J_{\bm a_{C\sm S}^{(S)}}\,:\,\RR^{|S|}\rightarrow \RR^{\lvert
    C\sm S\rvert \times |S|}$ denotes the Jacobian matrix of
  $\bm a_{C\sm S}^{(S)}$ evaluated at $\bm x_S$ and the square root in
  $\bm J_{\bm a_{C\sm S}^{(S)}}(\bm x_S)^{1/2}$ is interpreted
  element-wise.\ The limit vector $\bm Z_{V\sm \condv}^{(\condv)}$ in
  expression \eqref{eq:th11} is a zero-mean Gaussian graphical model
  with covariance matrix
  $\bm R_{V\sm v}^{(\condv)} = \left [ 2 \rho_{i\condv} \rho_{j\condv}
    ( \rho_{ij}-\rho_{i\condv}\rho_{j\condv} ) \right ]_{i,j \in V\sm
    \condv}$.
  The inverse of the covariance matrix is
  $\bm R_{V\sm v}^{(\condv)} = \text{diag}(1/(\sqrt{2}\bm
  \rho_{V\sm\condv,\condv})) \bm Q_{V\sm\condv} \text{diag}(1/
  (\sqrt{2} \bm \rho_{V\sm\condv,\condv}))$ where $\text{diag}$
  denotes the diagonal matrix, $\bm Q_{V\sm \condv}$ is the matrix
  $\bm Q$ without its $\condv$th row and $\condv$th column. 
  \color{black}
  
  Here we emphasize that convergence for $V \sm \condv$ can also be
  obtained directly using the results of \cite{hefftawn04} since for
  this graphical model, the joint distribution is specified directly
  on the entire graph $\GG$.\ This is different from the graphical
  model presented in Section~\ref{eq:HR_cliques_example} for which the
  joint distribution on the graph $\GG$ was specified indirectly via
  its lower-dimensional clique distributions.

  \subsection{Block graph with cliques of differing extremal dependence}
\label{sec:blockex}
A valid decomposable graphical model may be constructed with
individual clique distributions belonging to different families.  Such
a construction is possible if the distributions of the subvectors
corresponding to the separators are identical within each clique
distribution, that is
$F_{\bm X_C}(\bm x_S,\infty) = F_{\bm X_{C'}}(\bm x_S,\infty)$ for
$C \cap C' =S \neq \varnothing$ where $C,C' \in \mathcal{M}(\GG)$.\
This constraint is trivially satisfied for a block graph with
standardized univariate marginal distributions.\ Consider the subgraph
that is induced by the vertices $\{1,2,3,4,5\}$ in the graph presented
in Figure \eqref{fig:sepsimple}.\ This graph consists of the two
cliques $C_1 = \{1,2\}$, $C_2 = \{2,3,4,5\}$.\

Assume that $\bm X_{E, C_1}$ follows a bivariate max-stable \HR \
distribution with unit exponential marginal distributions and a
variogram matrix
$\bm \Gamma= (\Gamma_{ij})_{i,j=1}^2\in \RR^{2\times 2}$, with
$\Gamma_{12}>0$ and $\Gamma_{11}=\Gamma_{22}=0$. Assume further that
$\bm X_{E, C_2}$ follows a Gaussian distribution, also with unit
exponential univariate marginals and covariance matrix given by
$[\rho_{i,j}]_{i,j \in \{2,\dots, 5\}}$.\ Let $\condv=1$.\ Our
analysis of the examples in Sections \ref{eq:HR_cliques_example} and
\ref{sec:gaussian}, gives that as $t \rightarrow \infty$,
\begin{IEEEeqnarray*}{lCl}
  \PP \left (X_{E,2} - X_{E,1}  \leq z^{(1)}_2\ \big | \ X_{E,1} =t \right ) \wk \Phi( z^{(1)}_2, \ -\Gamma_{12}, \ 2\Gamma_{12}),  \\
  \PP \left ( \frac{\bm X_{E,\{3,4,5\}} - \bm
      \rho^2_{2,\{3,4,5\}}X_{E,2}} {{X_{E,2}}^{1/2}} \leq \bm
    \epsilon^{(2)}_{345} \ \bigg | \ X_{E,2} =t+z^{(1)}_2 \right ) \wk
  \Phi ( \bm \epsilon^{(2)}_{345},0 , \bm \Sigma^{(2)}),
\end{IEEEeqnarray*}
where
$\bm \Sigma^{(2)} = [ 2 \rho_{2,i}\rho_{2,j}(\rho_{i,j}
-\rho_{2,i}\rho_{2,j} ) ]_{i,j \in \{3,4, 5\}}$.  The associated
random vector $\bm \epsilon^{(2)}_{\{3,4,5\}}$ and the random variable
$Z^{(1)}_2$ are independent of each other.\ Assumption $\AB$ is
satisfied if we can find functions $\bm a^{(1)}_{\{3,4,5\}}(t)$ and
$\bm b^{(1)}_{\{3,4,5\}}(t)$ such that
\begin{IEEEeqnarray}{lCl}
  \lim_{\thres \to \infty} \frac{(\thres + z_2^{(1)})\bm
    \rho^2_{2,\{3,4,5\}} - \bm a^{(1)}_{\{3,4,5\}}(t)}{\bm
    b^{(1)}_{\{3,4,5\}}(t)} \in \RR^3 \quad \text{and} \quad
  \lim_{\thres \to \infty} \frac{{(\thres + {z_2^{(1)})^{1/2} }} \bm
    1_3}{\bm b^{(1)}_{\{3,4,5\}}(t)} > 0.
\label{eq:mconv}
\end{IEEEeqnarray}
The existence of the limits in expression \eqref{eq:mconv} is assured
if we let $\bm b^{(1)}_{345}(t) = {t}^{1/2} \bm 1_3$ and
$\bm a^{(1)}_{345}(t) = t \bm \rho^2_{2,\{3,4,5\}}$.\ In this case we
obtain, as $\thres \rightarrow \infty$, a non-degenerate Gaussian
distribution $G^{(1)}_{\{2,3,4,5\}}$ with mean $(-\Gamma_{12}, 0,0,0)$
and covariance consisting of a diagonal block matrix with elements
$ 2\Gamma_{12}$ and $\bm \Sigma^{(2)}$.

A similar construction with $\condv \in C_2 \sm \{2\}$ fails since the
convergence required in Assumption $\AB$ and illustrated in
\eqref{eq:mconv} cannot be ensured by any choice of
$\bm a^{(\condv)}_{C_1}(t)$ and $\bm b^{(\condv)}_{C_1}(t)$.\ The
reason for this is the different scale normalization required in the
two clique distributions.\ To illustrate this problem, let $\condv =3$
and consider the transition from clique 2 to clique 1.\ The Assumption
$\AB$ requires that
\begin{IEEEeqnarray*}{lCl}
\lim_{\thres \to \infty}  \frac{\rho^2_{2,3}\thres +  {t}^{1/2}z_3^{(2)}  -  a^{(3)}_{1}(t)}{ b^{(3)}_{1}(t)}  \in \RR \quad \text{and} \quad 
\lim_{\thres \to \infty} \frac{1}{b^3_{1}(t)} > 0.
\end{IEEEeqnarray*}
The second relationship is satisfied if $b^3_{1}(t) = \beta$ for all
$t>0$ where $\beta$ is a positive constant, but if this is the case
then the there is no solution for $ a^{(3)}_{1}(t)$.\ However, in this
case, we can still apply Theorem \ref{prop:block} and find that as
$\thres \rightarrow \infty$,
\begin{IEEEeqnarray*}{lCl}
 \PP  \left ( X_{E,1} - X_{E,2} \leq \epsilon^2_1, \frac{\bm X_{E,\{2,4,5\}} - \bm \rho^2_{3,\{2,4,5\}}X_{E,3}} {\sqrt{X_{E,3} }  }  \leq \bm z^{(3)}_{2,4,5}\ \bigg | \ X_{E,3} =t  \right ) & \wk   \Phi \left ( \{ \epsilon^{(2)}_1, \bm z^{(3)}_{2,4,5} \}  ,  \bm m, \bm S \right  ). \\
  \label{eq:randomnorm}
\end{IEEEeqnarray*}
The random elements $\varepsilon^{(2)}_1$ and
$\bm Z^{(3)}_{\{2,4,5\}}$ are independent of each other and follow a
joint Gaussian distribution with parameters
$\bm m = (-\Gamma_{12}, 0,0,0)$ and a covariance matrix $\bm S$
consisting of a diagonal block matrix with elements $ 2\Gamma_{12}$
and $\bm \Sigma^{(3)}$, where
$\Sigma^{(3)} = [ 2 \rho_{3,i}\rho_{3,j}(\rho_{i,j}
-\rho_{3,i}\rho_{3,j} ) ]_{i,j \in \{2,4, 5\}}$.
\begin{appendices}
\section{Graph theory and notation}
The material in this section is taken from \cite{grim18}.\ A graph
$\GG$ is a pair $(V, E)$, where $V= \{v_1, \dots, v_d\}$ is a finite
set of elements termed vertices, and $E\subseteq V\times V$.\ Each
element of $E$ is an unordered pair of vertices $v_i, v_j$ denoted by
$(v_i, v_j)$.\ Two edges with the same vertex-pairs are said to be in
parallel, and edges of the form $(v_i, v_i)$ are called loops.\ The
graphs in this paper contain neither parallel edges, nor loops and are
all assumed undirected since the edges are unordered pairs.\ Two
vertices $v_i$ and $v_j$ in $V$ are said to be adjacent if
$e=(v_i,v_j) \in E $.\ In this case, $v_i$ and $v_j$ are termed the
end vertices of $e$.\ A path of $G$ is defined as an alternating
sequence $v_0, e_0, v_1, e_2,\dots, e_{n-1},v_n$ of distinct vertices
$v_i$ and edges $e_i=(v_i, v_{i+1})$.\ Such a path has length $n$ and
it is said to connect $v_0$ and $v_{n}$.\ A path of length $0$ is
called a trivial path.\
A cycle of $\GG$ is an alternating sequence
$v_0, e_0, v_1, \dots, e_{n-1}, v_n, e_n, v_0$ such that
$v_0, e_0, v_1, \dots, e_{n-1}, v_n$ is a path and
$e_n=(v_{n}, v_{0})$.\ The graph-theoretic distance from $v$ to $v'$
is defined to be the number of edges in a shortest path of $\GG$ from
$v$ to $v'$.\ We write $v \leftrightsquigarrow v'$ if there exists a
path connecting $v$ and $v'$. The relation $\leftrightsquigarrow$ is
an equivalence relation, and its equivalence classes are called
components of $\GG$. The graph $\GG$ is termed connected if it has a
unique component and tree if in addition it contains no cycle. A graph
is termed chordal if every cycle of length greater than four has an
edge between two vertices that are non-consecutive in the cyle.\ A
subgraph of a graph $\GG$ is a graph $\mathcal{H} = (W, F)$ with
$W\subseteq V$ and $F\subseteq E\cap (W\times
W)\color{black}$.\ 
The subgraph $\mathcal{H}$ is a spanning tree of $\GG$ if $V=W$ and
$\mathcal{H}$ is a tree.\ Any connected graph $\GG$ has a spanning
tree.\ A vertex-induced subgraph $\GG_W$ is $(W, E\cap (W\times W))$.\
A subgraph is said to be complete if all its vertices are adjacent.\ A
complete subgraph of $\GG$ is called a clique, and a clique is maximal
if no strict superset is a clique.
 \label{app:graph}
  \section{Proofs}
  \label{app:proofs}
  \renewcommand{\theequation}{A.\arabic{equation}}
  \setcounter{equation}{0}

  \subsection{Proof of Proposition 1}
  \label{app:proof_prop_1}
  \begin{proof}
    Since $\bm X$ admits a Lebesgue density and its conditional
    independence graph is assumed decomposable, we have 
    \begin{IEEEeqnarray}{rCl}
      n^{d+1} f(n \bm x)&=& n^{d+1} f_{C_1}(n \bm x_{C_1})
      \,\prod_{i=2}^N \frac{
      f_{C_{\ttime}}(n \bm x_{C_i})}{f_{S_i}(n \bm x_{S_i})}
      =n^{d+1+K+|C_1|+1} f_{C_1}(n \bm x_{C_1})
      \,\prod_{i=2}^N \frac{n^{|C_i|+1} f_{C_i}(n \bm x_{C_i})
      }{n^{|S_i|+1} \,f_{S_i}(n \bm x_{S_i})},
      \label{eq:density_rescaled}
    \end{IEEEeqnarray}
    where $K=-1-|C_1|-\sum_{i=2}^N [|C_i| - |S_i|]$.\ Since each
    $\nu^{(C_i)}$ places all the mass in the interior of $\EE^{C_i}$
    and $\Lambda^{(C_i)}$ is assumed differentiable, the spectral
    measure $H^{(C_i)}$ that is associated with $\Lambda^{(C_i)}$
    admits a density $h_{C_i}$ in the interior of $\SSS_{{|C_i|-1}}$
    \citep{coles1991modelling}.\ A direct implementation of Theorem
    \ref{theorem:main} from the paper just cited gives
    \begin{equation}
      \lambda^{(C_i)}(\bm y_{C_i})=\frac{\partial^{|C_i|}}{\partial \bm y_{C_i}} \Lambda^{(C_i)}(\bm y_{C_i})=  \lVert  \bm
      y_{C_i}\rVert^{-(|C_i|+1)} h^{C_1}(\bm y_{C_i}/\lVert \bm y_{C_i}.
      \rVert)
      \label{eq:Lambda_deriv}
    \end{equation}
    Combining \eqref{eq:density_rescaled} with \eqref{eq:Lambda_deriv},
    and with the fact that for any decomposable graph $\GG$, the
    relation $|C_1| + \sum_{i=2}^N [|C_i| - |S_i|]=d$ always holds true,
    we get $n^{d+1}f_{\bm X}(n\,\bm x)\rightarrow \lambda(\bm x)$
    outside a set of Lebesgue measure zero, with $\lambda$ defined
    in expression \eqref{eq:lambda_dgm}.

    From Scheffe's Lemma \citep[][pg.\ 29]{bill08}, it follows that as
    $n\rightarrow\infty$ 
    \begin{IEEEeqnarray*}{rCl}
      n \,( 1-F_{\bm X}(n\,\bm y)) = \int_{[\bm 0, \bm y]^c}
      n^{d+1} f_{\bm X}(n\,\bm x)~\dd\bm x \rightarrow
      \int_{[\bm 0, \bm y]^c} \lambda(\bm x)\,\dd\bm
      x =: \Lambda(\bm y),
    \end{IEEEeqnarray*}
    at continuity points $\bm y \in \EE$ of the limit function
    $\Lambda$.\ Therefore, $F^n(n\bm x) \wk \exp(-\Lambda(\bm x))$.
  \end{proof}

  \subsection{Proof of Proposition 2}
  \label{app:proof_prop_2}

  \begin{proof}
    We have from the probability of union and intersection,
    \begin{IEEEeqnarray*}{rCl}
      \PP\left(\bm X_{C}/\level \in B\right) &=& \PP\big( \bm X_{C \sm S} / \level \in \EE^{C \sm S} \cap  \bm X_S / \level \in [0, \bm x_{S}]^c \big)
      = \PP \left (  \bm X_S / \level \in [0, \bm x_{S}]^c  \right )
      = \PP\left(\bm X_{C'}/\level \in B'\right), 
    \end{IEEEeqnarray*}
    since $\bm X_S$ has the same distribution as a subvector of either
    $\bm X_C $ or $\bm X_{C'}$.  Hence
    $\lim_{\thres \rightarrow \infty} \thres \PP\left(\bm X_{C}/\level
      \in B\right) = \lim_{\thres \rightarrow \infty} \thres
    \PP\left(\bm X_{C'}/\level \in B'\right)= \lim_{\thres \rightarrow
      \infty} \thres \PP \left ( \bm X_S / \level \in [0, \bm x_{S}]^c
    \right )$.
  \end{proof}
  \subsection{Subsidiary results}
  \begin{lemma}[\cite{kulik2015heavy}] \normalfont Let $\{\mu_n\}$ be a
    sequence of probability measures on $\RR^d$ which converges weakly
    to $\mu$ as $n \rightarrow \infty$.\ Let $\{\phi_n\}$ be a uniformly
    bounded sequence of measurable functions on $\RR^d$ such that
    $\phi_n$ converges to a continuous function $\phi$ where the
    convergence is uniform on compact sets of $\RR^d$.\ Then
    $\mu_n(\phi_n) \rightarrow \mu(\phi) $ as $n \rightarrow \infty$.
    \label{lemma:ks}
  \end{lemma}
  \begin{lemma}\normalfont
    Let $\mu_{\thres, \bm x}$ be a family of measures on $\RR^d$, where
    $\thres \in \RR$ and $\bm x \in \RR^{{\dims}}$.\ Assume that $\mu_{\thres, \bm x}$
    on $\RR^{d}$ converges weakly to the measure $\mu$ where the
    convergence is uniform on compact sets of $\bm x \in \RR^{{\dims}}$.\
    Then given a compact set, $C \subset \RR^{{\dims}}$, for any $\delta >0$,
    there exists a compact set $C_{\delta} \subset \RR^d$ such that for
    all $\thres>U_{\delta}$ and for all $\bm x \in C$ we have
    $\mu_{\thres, \bm x}(C_{\delta})> 1- \delta$.
    \label{lemma:tight}
  \end{lemma}
  \begin{proof}
    Suppose that the statement is not true, then for all compact sets
    $C \subset \RR^d$ we can choose an arbitrarily large $\level$ and an
    element $\bm x \in C$ such that
    $\mu_{\thres, \bm x}(\chi) \leq 1- \delta$.\ Choose a sequence
    $\thres_n \in \RR$ with $\thres_n>n$ and another sequence $\bm x_n\in C$ such
    that $\mu_{\thres_n, \bm x_n}(C) \leq 1- \delta$ for some $\delta$.\
    Since  $\thres_n \rightarrow \infty$ as $n \rightarrow \infty$, and since
     $\mu_{\thres, \bm x}(d\bm y) \xrightarrow{w} \mu(d \bm y)$ uniformly in the compact set $C$, the sequence of measures $\mu_{\thres_n, \bm x_n}$ must
    converge weakly to $\mu$.\  But if this sequence converges weakly
    then there exists some $C_{\delta}$ and an integer $N$ such that
    $\mu_{\thres_n, \bm x_n}(\chi_{\delta}) > 1-\delta $ for all $n>N$
    \citep[p. 105]{alma9924545702502466}, which is a contradiction.
  \end{proof}
  \begin{lemma}[Partition of unity, \cite{alma99309453502466} p. 251]
    \normalfont Let $C$ be a compact set in $\RR^d$ with a finite open cover
    $\{U_i\}_{i \in 1, \dots, I}$,  then there exists a set of continuous functions
    $\{k_j(\bm x)\}_{j=1,\dots,J}$, $\bm x \in C$, such that,
    \begin{enumerate}[wide=0\parindent]
    \item[$(i)$] for each $k_j(\bm x)$, $0 \leq k_j(\bm x) \leq 1$ for all
      $\bm x \in C$;
    \item[$(ii)$] for each function
      $k_j$ there is a cover set $U_{i_j}$ such that  $k_j(\bm x)=0$ when $\bm x \notin U_{i_j}$;
    \item[$(iii)$] $\sum_{j=1}^J k_j(\bm x) =1$, for all $\bm x \in C$.
    \end{enumerate}
    \label{lemma:partition}
  \end{lemma}
  \begin{lemma} \normalfont Let $\mu_{\thres, \bm x}$ be a family of
    measures on $\RR^d$, where $\thres \in \RR$ and
    $\bm x \in \RR^{{\dims}}$.\ Assume that $\mu_{\thres, \bm x} \wk \mu$
    uniformly in compact sets of $\bm x \in \RR^{{\dims}}$ as
    $\thres \rightarrow \infty$.\ Let
    $\bm \psi: \RR^{{\dims}} \rightarrow \RR^{d}$ and
    $\phi: \RR^{{\dims}} \rightarrow \RR^{d}$ be continuous functions such
    that $\phi(\bm x)>0$ for all $\bm x \in \RR^{{\dims}}$.\ Then, for any
    $f \in C_b(\RR^d)$
    \begin{equation*}
      \int_{\RR^{d}} f(\bm \psi(\bm x) + \phi(\bm x) \bm y) \mu_{\thres, \bm x}(d \bm y) \rightarrow \int_{\RR^{d}} f(\bm \psi(\bm x) + \phi(\bm x)\bm y) \mu(d \bm y), \quad \text{as } \thres \rightarrow \infty,
    \end{equation*}
    and the convergence is uniform on
    compact sets in the variable $\bm x \in \RR^{{\dims}}$.
    \label{lemma:secondterm}
  \end{lemma}
  \begin{proof}
    It suffices to show that given any $\epsilon > 0$, there exists $U_\epsilon$ such that
    \begin{equation}
      \sup_{\bm x_1 \in \mathcal{X}_1} \sup_{\bm x_2 \in \mathcal{X}_2}
      \left| \int_{\RR^{{\dims}}} f(\psi(\bm x_1) + \phi(\bm x_1) \bm y)
        \mu_{\thres, \bm x_2}(d \bm y) - \int_{\RR^{{\dims}}} f(\bm \psi(\bm x_1) +
        \phi(\bm x_1) \bm y) \mu(d \bm y) \right | < \epsilon,
      \label{eq:secondterm}
    \end{equation}
    for all $\thres>U_\epsilon$, where the sets
    $\mathcal{X}_1 \subset \mathcal{X}_2 \subset \RR^{{\dims}}$ are compact.\
    From Lemma \ref{lemma:tight}, there exists a compact set
    $\mathcal{Y} \subset \RR^d $ such that
    $\mu_{\thres, \bm x}(\mathcal{Y})>1-\epsilon$ for all $U>U_1$ for some
    $U_1$.\ The first term in expression \eqref{eq:secondterm} equals,
    \begin{equation*}
      \int_{\mathcal{Y}} f(\bm \psi(\bm x_1) + \phi(\bm x_1) \bm y) \mu_{\thres, \bm x_2}(d \bm y)  + \int_{\RR^{{\dims}} \sm \mathcal{Y} } f(\bm \psi(\bm x_1) + \phi(\bm x_1) \bm y) \mu_{\thres, \bm x_2}(d \bm y).
    \end{equation*}
    Consider the function
    $h(\bm x_1, \bm y) = \bm \psi(\bm x_1) + \phi(\bm x_1) \bm y $.\ Since
    both $ \bm x_1 \in \mathcal{X}_1 $ and $\bm y \in \mathcal{Y}$ are
    restricted to compact sets and since $\bm \psi$ and $\phi$ are both
    continuous, this function has a compact image
    $h(\mathcal{X}_1, \mathcal{Y}) \subset \RR^d$.\ The restriction of
    the continuous function $f$ to $h(\mathcal{X}_1, \mathcal{Y})$ is
    therefore uniformly continuous.\ Hence for any $\epsilon>0$, there
    exists a $\delta>0$ such that for all
    $\xi, \xi' \in h(\mathcal{X}_1, \mathcal{Y})$ where
    $|\xi - \xi|<\delta$ we have $|f(\xi) - f(\xi')|< \epsilon$.  In addition, function $f(h(\bm x_1, \bm y))$  is  equicontinuous in the variable $\bm y$ on the compact set $\mathcal{X}_1$. \ \\ Let
    $\sup_{\bm x_1 \in \mathcal{X}_1}|\bm \phi(\bm x_1)| = M_{\phi}$,
    then,
    \begin{equation}
      |h(\bm x_1, \bm y)-h(\bm x_1, \bm y')| = |\bm \phi(\bm x_1) (\bm y - \bm y')| \leq M_{\phi} |\bm y - \bm y'|.
      \label{eq:equi}
    \end{equation}
    Consequently, from the uniform continuity of $f$ on $h(\mathcal{X}_1, \mathcal{Y})$,  for all $\bm y$ and $\bm y'$ such that
    $|\bm y - \bm y'|< \delta/M_{\phi}$, we have
    $|f(h(\bm x_1, \bm y))-f(h(\bm x_1, \bm y'))|<\epsilon$ for all
    $\bm x_1 \in \mathcal{X}_1$.\ The set $\mathcal{Y} \subset \RR^d$ is
    compact so we can find a finite cover consisting of open spheres of
    radius $\delta/M_{\phi}$ with centres $\{\bm y^{(i)} \}$,
    $i=1, \dots, I$, such that for any point $\bm y$ within such a sphere we have $|f(h(\bm x_1, \bm y))-f(h(\bm x_1, \bm y^{(i)}))| <\epsilon$.\ So
    for all of the sphere centres $\bm y^{(i)}$,
    $f(h(\bm x_1, \bm y))= f\left (h(\bm x_1, \bm y^{(i)} \right) + \left [ f(h(\bm
      x_1, \bm y)) - f(h(\bm x_1, \bm y^{(i)}) \right ] = f\left (h(\bm x_1,
    \bm y^{(i)} \right ) + R(\bm x_1, \bm y, \bm y^{(i)})$, where by
    construction $|R(\bm x_1, \bm y, \bm y')| < \epsilon$ for all
    $\bm x_1$.\  By Lemma \ref{lemma:partition} there exists a (finite)
    partition of unity on this cover $\{k_j(\bm y)\}_{j=1,\dots,J}$,
    such that $\sum_{j=1}^J k_j(\bm y)=1$ for all
    $\bm y \in \mathcal{Y}$ with the support of each $k_j(\bm y)$ within
    a cover sphere centred on $\bm y^{(j)}$ so that for all
    $\bm y \in \mathcal{Y}$,
    \begin{equation}
      f(h(\bm x_1, \bm y)) =  f(h(\bm x_1, \bm y))\sum_{j=1}^J k_j(\bm y)
      = \sum_{j=1}^J \left [f(h(\bm x_1, \bm y^{(j)})) + R(\bm x_1, \bm y, \bm y^{(j)} ) \right] k_j(\bm y),
      \label{eq:part2}
    \end{equation}
    since $k_j(\bm y)=0$  outside the open cover sphere centered on $\bm y^{(j)}$.  Returning to expression \eqref{eq:secondterm} and applying the
    triangle inequality,
    \begin{IEEEeqnarray}{lCl}
      & &\hspace{-40pt}\left| \int_{\RR^d} f(\bm \psi(\bm x_1) + \phi(\bm x_1) \bm y)
        \mu_{\level, \bm x_2}\color{black}(d \bm y) -
        \int_{\RR^d} f(\bm \psi(\bm x_1) + \phi(\bm x_1) \bm y) \mu(d \bm y)
      \right | \leq \nonumber\\ &&  \hspace{-20pt}\leq \left| \int_{\mathcal{Y}}
        f(\bm \psi(\bm x_1) + \phi(\bm x_1) \bm y) \mu_{\thres, \bm x_2}(d
        \bm y) - \int_{\mathcal{Y}} f(\bm \psi(\bm x_1) + \phi(\bm x_1) \bm
        y) \mu(d \bm y) \right | + \nonumber\\ &+&  \left|
        \int_{\RR^d \sm \mathcal{Y}} f(\bm \psi(\bm x_1) + \phi(\bm x_1)
        \bm y) \mu_{\thres, \bm x_2}(d \bm y) - \int_{\RR^d \sm
          \mathcal{Y}} f(\bm \psi(\bm x_1) + \phi(\bm x_1) \bm y) \mu(d \bm
        y) \right |.
      \label{eq:secondterm1}
    \end{IEEEeqnarray}
    \color{black} Let $M_f$ be the (global) supremum of the function
    $f \in C_b(\RR^d)$.\ From the construction in equation
    \eqref{eq:part2} we can place bounds on the first term of the
    right-hand side of expression \eqref{eq:secondterm1} so that,
    \begin{IEEEeqnarray*}{lCl}
      &\sup_{\bm x_1 \in \mathcal{X}_1} & \sup_{\bm x_2 \in
        \mathcal{X}_2}\,\left| \int_{\mathcal{Y}} f(\bm \psi(\bm x_1) +
        \phi(\bm x_1) \bm y) \mu_{\thres, \bm x_2}(d \bm y) -
        \int_{\mathcal{Y}} f(\bm \psi(\bm x_1) +
        \phi(\bm x_1) \bm y) \mu(d \bm y) \,\right| \leq \\\\
      &\leq& M_f \sup_{\bm x_2 \in \mathcal{X}_2} \sum_{j=1}^J \,\left|
        \int_{\mathcal{Y}} k_j(\bm y)\mu_{\thres, \bm x_2}(d \bm y) -
        \int_{\mathcal{Y}} k_j(\bm y)\mu(d \bm y) \,\right| \quad + \\\\
      &+& \sum_{j=1}^J \sup_{\bm x_1 \in \mathcal{X}_1} \sup_{\bm x_2
        \in \mathcal{X}_2}\,\left| \int_{\mathcal{Y}} R(\bm x_1, \bm y,
        \bm y^{(j)} ) \mu_{\level, \bm x_2}\color{black}(d \bm
        y) \,\right| + \sum_{j=1}^J \sup_{\bm x_1 \in \mathcal{X}_1}
      \sup_{\bm x_2 \in \mathcal{X}_2}\,\left| \int_{\mathcal{Y}} R(\bm
        x_1, \bm y, \bm y^{(j)} )\mu(d \bm y) \,\right|.
    \end{IEEEeqnarray*}
    Since $\mu_{\thres, \bm x_2}(d \bm y)$ converges weakly to
    $\mu(d \bm y)$, uniformly over compact sets of $\bm x_2$, and since
    the remainder term $R(\bm x_1, \bm y, \bm y^{(j)})$ satisfies
    $\left | R(\bm x_1, \bm y, \bm y^{(j)}) \right | < \epsilon$ on
    spheres centred at $\bm y^{(j)}$ of radius $\delta/M_{\phi}$, we can
    find $V_2$ such that for any $v> V_2$,
    \begin{equation*}
      \sup_{\bm x_1 \in \mathcal{X}_1} \sup_{\bm x_2 \in \mathcal{X}_2}
      \left| \int_{\mathcal{Y}} f(\bm \psi(\bm x_1) + \phi(\bm x_1) \bm y)
        \mu_{\thres, \bm x_2}(d \bm y) - \int_{\mathcal{Y}} f(\bm \psi(\bm x_1) +
        \phi(\bm x_1) \bm y) \mu(d \bm y) \right | < J M_f \epsilon +2J  \epsilon.
    \end{equation*}
    Also, using the definition of the compact set $\mathcal{Y}$ and by
    another application of the triangle inequality, the second term on the
    right-hand side of expression \eqref{eq:secondterm1} satisfies, for $\thres>U_1$,
    \begin{equation*}
      \sup_{\bm x_1 \in \mathcal{X}_1} \sup_{\bm x_2 \in \mathcal{X}_2}
      \left| \int_{\RR^d \sm \mathcal{Y}} f(\bm \psi(\bm x_1) + \phi(\bm x_1) \bm y)
        \mu_{\thres, \bm x_2}(d \bm y) - \int_{\RR^d \sm \mathcal{Y}} f(\bm \psi(\bm x_1) +
        \phi(\bm x_1) \bm y) \mu(d \bm y) \right | < 2M_f \epsilon.
    \end{equation*}
    Hence by setting $U_\epsilon > \max(U_1,U_2)$ we obtain an arbitrarily small bound on expression \eqref{eq:secondterm}.
  \end{proof}
  \begin{lemma} \normalfont Let $C$ and $S$ denote a clique and a
    separating set, respectively, of a decomposable graphical model.\
    Let $\bm X_E$ be a random vector that satisfies Assumptions
    $\AAA$, $\AA$ and $\AB$. Let the functionals
    $\bm a^{(\condv)}_{S}(\level)$ and $\bm b^{(\condv)}_{S}(\level)$
    be as defined in $\AB$.\ Then, for every
    $f \in C_b(\RR^{|C \sm S|}) $,
    \begin{multline}
      \int_{ \RR^{|C \sm S|}} f(\bm y) \PP\left( \frac{ \bm X_{E,C \sm S}- \bm a^{(\condv)}_{C \sm S}(\level)}{\bm b^{(\condv)}_{C \sm S}(\level)} \in d\bm y \Big | \frac{\bm  X_{E,S} - \bm a^{(\condv)}_{S}(\level)} {\bm b^{(\condv)}_{S}(\level)} = \bm x  \right) \rightarrow \\
      \int_{ \RR^{|C \sm S|}} f \left ( \bm  \psi^{(S)}_{C \sm S}(\bm x)   + \bm \phi^{(S)}_{C \sm S}(\bm x) \bm y  \right ) G^{(S)}_{C \sm S}(d\bm y),
      \text{ as } \level \to \infty,
      \label{eq:lemmawc}
    \end{multline}
    and the convergence is uniform on compact sets of the variable
    $\bm x \in \RR^{ |S|}$.
    \label{lemma:01}
  \end{lemma}
  \begin{proof}
    From \eqref{eq:remainder}, we have
    $\bm a^{(\condv)}_{C \sm S}(\level) + \bm b^{(\condv)}_{C \sm
      S}(\level)\, \bm y$ equals
    \begin{equation*}
      \bm a^{(S)}_{C \sm S}\left (\bm T_S^{(\condv)}(\bm x, \level) \right) +
      \bm b^{(S)}_{C \sm S}\left (\bm T_S^{(\condv)}(\bm x, \level) \right)
      \left [ \bm A^{(S)}_{C \sm S}(\bm x, \level) +
        \left( \bm 1- \bm  B^{(S)}_{C \sm S} (\bm x, \level)\right)
        \frac{\bm y -  \bm \psi^{(S)}_{C \sm S} (\bm x)}
        {\bm \phi^{(S)}_{C \sm S} (\bm x)}\right].
    \end{equation*}
    Using this relationship, we perform a change of variable in
    expression \eqref{eq:lemmawc} to remove the dependence on
    $\bm a^{(\condv)}_{C \sm S}(\level)$ and
    $\bm b^{(\condv)}_{C \sm S}(\level)$, that is,
    \begin{IEEEeqnarray*}{rCl}
      &\int_{ \RR^{|C \sm S|}}& f \left ( \bm y \right)\PP\left(
        \frac{ \bm X_{C \sm S}- \bm a^{(\condv)}_{C \sm S}(\bm
          X_S)}{\bm b^{(\condv)}_{C \sm S}(\bm X_S)} \in d\bm y
        ~\Big|~ \bm X_{S} = \bm T_S^{(\condv)}(\bm x, \level)\right) =
      \int_{ \RR^{|C \sm S|}} f(\bm g(\bm x,\level) + \bm y \bm h(\bm
      x,\level)) 
      \mu_{\level, \bm x} (d \bm y) 
    \end{IEEEeqnarray*}
    where
    \begin{IEEEeqnarray*}{rCl}
      \mu_{\level, \bm x} (d \bm y) &=& \PP \left( \frac{ \bm X_{C \sm
            S}- \bm a^{(S)}_{C \sm S}\left (\bm T_S^{(\condv)}(\bm x, \level)
          \right )}{\bm b^{(S)}_{C \sm S}\left (T_S^{(\condv)}(\bm x,
            \level)\right )} \in d\bm y \Big | \bm X_{S} = \bm T_S^{(\condv)}(\bm x,
        \level) \right),\\
      \bm g(\bm x,\level) &=&   \bm \psi^{(S)}_{C \sm S}(\bm x) -\frac{ \bm \phi^{(S)}_{C \sm S}(\bm x) \bm  A^{(S)}_{C \sm S}(\bm x, \level)}{1- \bm  B^{(S)}_{C \sm S}(\bm x, \level)},\quad \text{and} \quad
      \bm h(\bm x,\level) = \frac{\bm \phi^{(S)}_{C \sm S}(\bm x) }{1- \bm  B^{(S)}_{C \sm S}(\bm x, \level)},
    \end{IEEEeqnarray*}
    Suppressing subscripts and superscripts for simplicity, we require
    that
    \begin{equation*}
      \sup_{\bm x \in \mathcal{X} }\, \left|\,\int_{ \RR^{|C \sm S|}} f \left( \bm g(\bm x,\level) + \bm y\bm h(\bm x,\level)  \right ) \mu_{\level, \bm x} (d \bm y) - \int_{ \RR^{|C \sm S|}} f \left (\bm \psi(\bm x)   +\bm \phi(\bm x) \bm y  \right ) G(d\bm y) \,\right| \rightarrow 0,
      \text{ as }  \level\rightarrow \infty,
    \end{equation*}
    where the set $\mathcal{X} \subset \RR^{|S|}$ is compact.\ Applying
    the triangle inequality gives,
    \begin{multline}
      \left |\int_{ \RR^{|C \sm S|}} f \left ( \bm g(\bm x,\level) + \bm y\bm h(\bm x,\level)  \right ) \mu_{\level, \bm x} (d \bm y) -  \int_{ \RR^{|C \sm S}} f \left (\bm \psi(\bm x)   + \bm \phi(\bm x) \bm y  \right ) G(d\bm y) \right | \\
      \leq \int_{ \RR^{|C \sm S|}} \left | f \left (  \bm g(\bm x,\level) + \bm y\bm h(\bm x,\level)  \right ) - f \left ( \bm \psi(\bm x)   + \bm \phi(\bm x) \bm y  \right )\ \right|\  \mu_{\level, \bm x}  (d \bm y)\\
      +\left | \int_{ \RR^{|C \sm S|}} f \left (\bm \psi(\bm x) \bm  +\phi(\bm x) \bm y  \right ) \mu_{\level, \bm x} - \int_{ \RR^{|C \sm S}}  \left ( \bm \psi(\bm x)   + \bm \phi(\bm x) \bm y  \right ) G(d\bm y)  \right |.
      \label{eq:ineq}
    \end{multline}
    We consider each term in the right-hand side of inequality \eqref{eq:ineq}
    separately.
    \begin{description}[wide=0\parindent]
    \item[First term.] Let
      $\mathcal{X}_1 \subset \mathcal{X}_2 \subset \RR^{|S|}$, with
      $\mathcal{X}_1$ and $\mathcal{X}_2$ both compact.\ We consider
      the more general supremum,
      \begin{equation}
        \sup_{\bm x_2 \in \mathcal{X}_2 } \sup_{\bm x_1 \in \mathcal{X}_1 }  \int_{ \RR^{|C \sm S|}} \left | f \left (  \bm g(\bm x_1,\level) + \bm y\bm h(\bm x_1,\level)  \right ) - f \left ( \bm \psi(\bm x_1)   + \bm \phi(\bm x_1) \bm y  \right )\ \right|\  \mu_{\level, \bm x_2}  (d \bm y)
        \label{eq:anotherfirst}
      \end{equation}
      From \color{black} Lemma \ref{lemma:ks}\color{black}, it follows
      that for every $\epsilon>0$ there exists a compact set
      $\mathcal{Y}$ and $U_1>0$ such that
      $\mu_{\level, \bm x_2}(\mathcal{Y}) > 1- \epsilon$ for
      $\level>U_1$.\ From Assumption $\AB$,
      $\bm g(\bm x_1,\level) \rightarrow \bm \psi(\bm x_1) $ and
      $\bm h(\bm x_1,\level) \rightarrow \bm \phi(\bm x_1) $,
      uniformly in compact sets in the variable
      $\bm x_1 \in \RR^{|S|}$.\ For $\bm y$ in the compact set
      $ \mathcal{Y}$,\ the continuity of $f$ on $\RR^{|C \sm S|}$
      implies uniform convergence of the terms involving the function
      $f \in C_b(\RR^{|C \sm S |})$, so that,
      \begin{equation*}
        \sup_{ \bm x_1 \in \mathcal{X}}
        \left | f  \left (  \bm g(\bm x_1,\level) + \bm y\bm h(\bm x_1,\level)   \right )
          -  f\left ( \bm \psi(\bm x_1) \bm  +\phi(\bm x_1) \bm y  \right ) \right |
        < \epsilon \quad \text{for all } \bm y \in \mathcal{Y},
        \quad \text{for all } t>U_2.  \quad
      \end{equation*}
      In addition, if we set
      $M_f:=\sup\{\,\left | f(\bm x)\right |\, :\,{\bm x \in \RR^{|C
          \sm S|}}\}$, then from the triangle inequality it follows
      that
      \begin{equation*}
        \Big | f  \left (   \bm g(\bm x_1,\level) + \bm y\bm h(\bm x_1,\level)   \right ) -  f \left (\bm \psi(\bm x_1) \bm  +\phi(\bm x_1) \bm y  \right ) \Big | <2M_f.
      \end{equation*}
      Hence, the supremum of the first term on the right-hand side of
      expression \eqref{eq:ineq} satisfies
      \begin{multline*}
        \sup_{ \bm x_2 \in \mathcal{X}_2}  \int_{ \RR^{|C \sm S|}} \left ( \sup_{ \bm x_1 \in \mathcal{X}_1} \left |f \left (  \bm g(\bm x_1,\level) + \bm y\bm h(\bm x_1,\level)  \right ) -  f \left (\bm \psi(\bm x_1) \bm  +\phi(\bm x_1) \bm y  \right ) \right | \right ) \mu_{\level, \bm x_2}  (d \bm y) \\
        \leq \sup_{ \bm x_2 \in \mathcal{X}_2}  \int_{ \mathcal{Y}}  \left ( \sup_{ \bm x_1 \in \mathcal{X}_1} \left |f \left (  \bm g(\bm x_1,\level) + \bm y\bm h(\bm x_1,\level)  \right ) -  f \left (\bm \psi(\bm x_1) \bm  +\phi(\bm x_1) \bm y  \right ) \right | \right ) \mu_{\level, \bm x_2}  (d \bm y) \\
        +\sup_{ \bm x_2 \in \mathcal{X}_2}  \int_{ \RR^{|C \sm S|} \sm \mathcal{Y}} \left ( \sup_{ \bm x_1 \in \mathcal{X}_1} \left |f \left (  \bm g(\bm x_1,\level) + \bm y\bm h(\bm x_1,\level)  \right ) -  f \left (\bm \psi(\bm x_1) \bm  +\phi(\bm x_1) \bm y  \right ) \right | \right ) \mu_{\level, \bm x_2}  (d \bm y)\\
        < \epsilon   \int_{\mathcal{Y}} \mu_{\level, \bm x_2}  (d \bm y)  +  2 M_f \epsilon <  (1+2M_f) \epsilon.
      \end{multline*}
      Hence, choosing $\level>\max(U_1, U_2)$, then the first term of
      expression \eqref{eq:ineq} can be made arbitrarily small.
    \item[Second term.] From Lemma \ref{lemma:secondterm}, we have
      that as $\level \rightarrow \infty$,
      \begin{equation*}
        \int_{ \RR^{|C \sm S|}} f \left (\bm \psi(\bm x) \bm  +\phi(\bm x) \bm y  \right ) \mu_{\level, \bm x} \rightarrow \int_{ \RR^{|C \sm S}} f \left ( \bm \psi(\bm x)   + \bm \phi(\bm x) \bm y  \right ) G(d\bm y),
      \end{equation*}
      uniformly on compact sets in the variable
      $\bm x \in \RR^{|S|}$.\ Hence, there exists $U_3>0$ such that,
      for all $\level>U_3$,
      \begin{equation*}
        \sup_{ \bm x \in \mathcal{X}} \left| \int_{ \RR^{|C \sm S|}} f \left (\bm \psi(\bm x) \bm  +\phi(\bm x) \bm y  \right ) \mu_{\level, \bm x} - \int_{ \RR^{|C \sm S}}  \left ( \bm \psi(\bm x)   + \bm \phi(\bm x) \bm y  \right ) G(d\bm y) \right| < \epsilon
        \label{eq:state1}
      \end{equation*}
    \end{description}
    Hence, for all $\level>\max(U_1, U_2,U_3)$, the supremum over
    $\bm x \in \mathcal{X}$ of the right-hand side of \eqref{eq:ineq}
    is bounded by $\epsilon + (1+2M_f) \epsilon = 2(1+M_f) \epsilon$,
    which proves the required result.
  \end{proof}
  \begin{remark} \normalfont The results in Lemma \ref{lemma:01} are
    equally valid for the real integral of a complex function $f$,
    with the modulus in the complex plane replacing the absolute value
    in the various inequalities.\ Similarly Lemma \ref{lemma:ks} is
    valid for complex functions by linearity of the (real) integral.
    \label{rem:complex}
  \end{remark}

  \subsection{Proof of Theorem 1}
  
  Let the graph $\GG= (V,E)$ have cliques $C_1,\dots,C_N$, which have
  been ordered according to the running intersection~\eqref{eq:ri}.\
  Without loss of generality, we assume that $\condv \in C_1$, see
  Remark~\ref{rem:rip}.\ Consider the measures,
  \begin{IEEEeqnarray*}{rCl}
    \mu_{\GG}^{\thres}(d\bm {z} )& = & \prod_{\ttime=2}^{N} \PP \left(
      \frac { \bm X_{C_{\ttime} \sm S_{\ttime}}- \bm
        a^{(\condv)}_{C_{\ttime} \sm S_{\ttime}}(X_{\condv})}{\bm
        b^{(\condv)}_{C_{\ttime} \sm S_{\ttime}}(X_{\condv})} \in d\bm
      z_{C_{\ttime} \sm S_{\ttime}} \ \Big | \ \frac{\bm
        X_{S_{\ttime}} - \bm a^{(\condv)}_{S_{\ttime}}(X_{\condv})}
      {\bm b^{(\condv)}_{S_{\ttime}}(X_{\condv})} = \bm z_{S_{\ttime}}
    \right) \\\\ & \times & \PP\left (\frac{\bm X_{C_1 \sm \condv}
        -\bm a^{(\condv)}_{C_1 \sm \condv } (X_{\condv})} {\bm
        b^{(\condv)}_{C_1 \sm \condv }(X_{\condv}) } \in d\bm z_{C_1
        \sm \condv } \ \bigg| \ X_{\condv} = \thres +
      \color{black}x_{\condv} \right )\, \PP\left
      (X_{\condv} - \level \in dx_{\condv} \ \bigg | \
      X_{\condv}>\level \right ),
  \end{IEEEeqnarray*}
  and 
  \begin{equation*}
    \mu_{\GG}(d\bm {z}) = \prod_{\ttime=1}^{N-1}
    G^{(S)}_{C_{\ttime}  \sm S_{\ttime}}
    \left ( \frac{d\bm z_{C_{\ttime}  \sm S_{\ttime}}  -
        \bm \psi^{(S_{\ttime})}_{C_{\ttime}  \sm S_{\ttime}}(\bm z_{S_{\ttime}})}
      {\phi^{(S_{\ttime})}_{C_{\ttime}  \sm S_{\ttime}}(\bm z_{S_{\ttime}})}
    \right)
    G^{(\condv)}_{C_1 \sm \condv}
    \left (d\bm z_{C_1 \sm \condv }\right ) \expdist(dx_{\condv}).
  \end{equation*}  
  We seek to show that $\mu_{\GG}^{\thres} \wk \mu_{\GG}$ as
  $\thres\rightarrow \infty$ by proving that the characteristic
  functions converge pointwise, that is,
  $\chi_{\mu_{\GG}^{\thres}}(\bm \omega) \rightarrow
  \chi_{\mu_{\GG}}(\bm \omega)$ as $\thres\rightarrow \infty$ for all
  $\bm \omega\in\RR^{d}$, where $\iota^2 = -1$,
  \begin{equation*}
    \chi_{\mu_{\GG}^{\thres}}(\bm \omega) = \int_{\RR^{d-1}\times [0,\infty)}
    \exp(\iota \, \bm \omega^\top \bm z)\mu_{\GG}^{\thres}(d\bm {z} ),\qquad \text{and}\qquad
    \chi_{\mu_{\GG}}(\bm \omega) = \int_{\RR^{d-1}\times [0,\infty)}
    \exp(\iota \, \bm \omega^\top  \bm z)\mu_{\GG}(d\bm {z} ),
  \end{equation*}
  We introduce the following shorthand notation for the measures
  \begin{IEEEeqnarray*}{rCl}
    \mu_{C_{\ttime},S_{\ttime},\condv}(d \bm z_{C_{\ttime} \sm
      S_{\ttime}}) & = & \PP \left( \frac{ \bm X_{C_{\ttime} \sm
          S_{\ttime}}- \bm a^{(\condv)}_{C_{\ttime} \sm
          S_{\ttime}}(X_{\condv})}{\bm b^{(\condv)}_{C_{\ttime} \sm
          S_{\ttime}}(X_{\condv})}
      \in d\bm z_{C_{\ttime}  \sm S_{\ttime}} \Big | \frac{\bm  X_{S_{\ttime}} - \bm a^{(\condv)}_{S_{\ttime}}(X_{\condv})} {\bm b^{(\condv)}_{S_{\ttime}}(X_{\condv})}  = \bm z_{S_{\ttime}}  \right) \quad \ttime = 2, \dots, N, \\\\
    \mu_{C_1 \sm \condv,\thres}(d\bm z_{C_1 \sm \condv}) &=&  \PP\left (\frac{\bm X_{C_1 \sm \condv} -\bm a^{(\condv)}_{C_1 \sm \condv } (X_{\condv})} {\bm b^{(\condv)}_{C_1 \sm \condv}(X_{\condv}) }  \in d\bm z_{C_1 \sm \condv} \ \bigg| \  X_{\condv} = \thres + x_{\condv} \right ), \qquad \text{and}\\\\
    \mu_{\condv,\thres}(dx_{\condv}) &=& \PP\left(X_{\condv} - \thres
      \in dx_{\condv} \ \bigg | \ X_{\condv}>\thres \right ).
  \end{IEEEeqnarray*}
  We proceed by induction on $\{1,\dots,N\}$,\ adding a clique to a
  graph to create a new graph with $N+1$ cliques.\ The hypothesis
  holds true for $N=1$, since such a graph consists of the single
  clique, $C_1$ and so the result follows from Assumption $\AA$.\
  
  We evaluate the characteristic function
  $\chi_{\mu_{\GG}^\thres}(\bm \omega)$ for the graph with $N+1$
  cliques through a nested integration containing the characteristic
  function for a graph with size $N$ cliques.  By conditional
  independence,
  \begin{multline}
    \chi_{\GG}^{\thres}(\omega) = \int_{
      \RR^{d-1} \times [0, \infty)}
    \exp ( \iota \bm \omega_{\GG }\cdot \bm z_{\GG } ) \\
    \times \left [ {\int_{ \RR^{|C \sm S|}}} \exp \left (\iota \bm t_{C \sm S} \cdot \bm z_{C_{N+1} \sm S_{N+1}} \right ) \mu_{C_{N+1},S_{N+1},\condv}(d \bm
      z_{C_{N+1} \sm S_{N+1}})\right ] \\
    \times \mu_{\GG \sm v ,\condv,\thres}(d\bm z_{\GG \sm v})
    \mu_{\condv,\thres}(dz_k).
    \label{eq:brackets}
  \end{multline}
  Concentrating on the terms contained within the square bracket, Lemma
  \ref{lemma:secondterm} and Remark \ref{rem:complex} imply that,
  \begin{equation*}
    {\int_{ \RR^{|C \sm S|}}} \exp \left (\iota \,\bm \omega_{C \sm S} \cdot  \bm z_{C \sm S} \right )  \mu_{C,S,k}(d \bm z_{C \sm S})  \to \int_{ \RR^{|C \sm S|}} \exp \left (\iota\, \bm \omega_{C  \sm S} \cdot \bm z_{C  \sm S} \right )  G^{(S)}_{C  \sm S}   \left (  \frac{d\bm z_{C  \sm S}  - \bm{\psi}^{(S)}_{C  \sm S}(\bm z_{S})  }{\bm{\phi}^{(S)}_{C  \sm S}(\bm z_{S}) } \right ),
  \end{equation*}
  as $\thres \rightarrow \infty$, where convergence is uniform on
  compact sets of $\RR^{|S|}$.\ Furthermore, since the term in the
  square brackets has real and imaginary parts that are bounded, we
  have by Lemma \ref{lemma:ks}, Remark \ref{rem:complex} and the
  inductive hypothesis that the right-hand side of equation
  \eqref{eq:brackets} tends to
  \begin{align*}
    \int_{ \RR^{d-1} \times [0, \infty)} \exp (\iota \bm \omega_{\GG} \cdot  \bm z_{\GG} )      \left [\int_{ \RR^{|s \sm c|}} \exp \left (\iota \bm \omega_{C  \sm S} \cdot  \bm z_{C  \sm S} \right )  G^{(S)}_{C  \sm S}   \left (  \frac{d\bm z_{C  \sm S}  - \bm{\psi}^{(S)}_{C  \sm S}(\bm {z_{S}})  }{\bm{\phi}^{(S)}_{C  \sm S}(\bm {z_{S}})}  \right ) \right ] 
    \mu_{\GG}(d \bm z_{\GG}).  \qed
  \end{align*}

  \subsection{Proof of Theorem \ref{prop:block}}
  We require first the following lemma.
  \begin{lemma} \normalfont Suppose that Assumption $\BA$ holds. Let
    $H^{1}(t) = t$ and for $\ttime=2,\dots,N$, let
    \[
      H^{\ttime}(t, y_{S_2}, \dots, y_{S_{\ttime}})=
      a^{S_{\ttime-1}}_{S_{\ttime}} \left \{ H^{\ttime-1}(t, y_{S_2},
        \dots, y_{S_{\ttime-1}} ) \right \}+
      b^{S_{\ttime-1}}_{S_{\ttime}} \left \{ H^{\ttime-1}(t, y_{S_2},
        \dots, y_{S_{\ttime-1}} ) \right \} y_{S_i} \in \RR.
    \]
    Then for any compact set $K^{\ttime} \subset \RR^{\ttime-1}$,
    $\ttime=1,\dots,N$, where $K^0=\emptyset$, and for arbitrarily
    large $T$, there exists $T_{\ttime}$ such that $t> T_{\ttime}$
    implies $H^{\ttime}(t, y_{S_2}, \dots, y_{S_{\ttime}}) > T$ for
    all $ (y_{S_2}, \dots, y_{S_{\ttime}}) \in K^{\ttime}$.
    \label{lemma:inf}
  \end{lemma}
  \begin{proof}[Proof of Lemma \ref{lemma:inf}]
    The proof is by induction.\ The claim is immediate for $\ttime=1$.
    Assume the claim holds true for $\ttime$ such that $1<\ttime<N $.\
    It suffices to prove the claim for the cylinder set
    $K^{\ttime+1}= K^{\ttime}\times [y_{\min}, y_{\max}]$ where
    $y_{S_{\ttime+1}} \in [y_{\min}, y_{\max}]$.\ We have
    $H^{\ttime+1}(t, y_{S_2}, \dots, y_{S_{\ttime+1}}) \geq
    a^{(S_{\ttime})}_{S_{\ttime+1}} \left \{ H^{\ttime}(t, y_{S_2},
      \dots, y_{S_{\ttime}} ) \right \}+
    b^{(S_{\ttime})}_{S_{\ttime+1}} \left \{ H^{\ttime}(t, y_{S_2},
      \dots, y_{S_{\ttime}} ) \right \} y_{\min}$ and by the induction
    hypothesis, $t > T_{\ttime}$ implies that
    $H^{\ttime}(t, y_{S_2}, \dots, y_{S_{\ttime}} ) >T$ where $T$ is
    arbitrarily large.\ From Assumption $\BA$ we know that
    $\bm a^{(S_{\ttime})}_{C_{\ttime+1} \sm S_{\ttime+1} }(\level)
    +a^{(S_{\ttime})}_{C_{\ttime+1} \sm S_{\ttime+1} }(\level) \bm
    y_{C_{\ttime+1} \sm S_{\ttime+1} } \rightarrow \infty \, \bm
    1_{|C_{\ttime}|-1} $ as $\level \rightarrow \infty$, hence we can
    choose $T_{\ttime+1}$ such that $t> T_{\ttime+1}$ implies
    $H^{\ttime+1}(t, y_{S_2}, \dots, y_{S_{\ttime+1}}) > T $.
  \end{proof}
  Consider the measures,
  \begin{IEEEeqnarray*}{rCl}
    \nu^{t,j}_{V }(d \bm y_{V \sm v}) &=& \PP \left ( \frac{ \bm
        X_{E,C_{1} \sm \condv } -\bm a^{(\condv)}_{C_{1} \sm \condv }
        (X_{E,\condv}) } {\bm b^{(\condv)}_{C_{2} \sm
          \condv}(X_{E,\condv})} \in d \bm y_{C_1 \sm \condv} \ \bigg
      | \ X_{E, \condv} = t \right) \nu_{C_{\ttime} \sm
      S_{\ttime}|S_{\ttime-1}}(d \bm y_{C_{\ttime} \sm
      S_{\ttime}}),
    \qquad \text{and} \\
    \nu^j_{V \sm \condv }(d \bm y_{ V \sm \condv }) &= &
    G^{(\condv)}_{C_1 \sm \condv}(d \bm y_{C_1 \sm \condv})
    \prod_{\ttime =2}^j G^{(S_{\ttime})}_{C_{\ttime} \sm S_{\ttime}}(d
    \bm y_{C_{\ttime} \sm S_{\ttime}}),
  \end{IEEEeqnarray*}
  where
  \[
    \nu_{C_{\ttime} \sm S_{\ttime}|S_{\ttime-1}}(d \bm y_{C_{\ttime}
      \sm S_{\ttime}}) = \PP \left ( \frac{ \bm X_{E,C_{\ttime} \sm
          S_{\ttime}}- \bm a^{(S_{\ttime})}_{C_{\ttime} \sm
          S_{\ttime}}( X_{E,S_{\ttime}})}{\bm
        b^{(S_{\ttime})}_{C_{\ttime} \sm S_{\ttime}}(
        X_{E,S_{\ttime}}) }\in d \bm y_{C_{\ttime} \sm S_{\ttime}} \
      ~\bigg|~\frac{X_{E, S_{\ttime}}-
        a^{S_{\ttime-1}}_{S_{\ttime}}(X_{E,S_{\ttime-1}})}{b^{S_{\ttime-1}}_{S_{\ttime}}(X_{E,S_{\ttime-1}})}
      = y_{ S_{\ttime}}\right ).
  \]
  Then we need to show that
  $ \nu^{t,j}_{V \setminus \condv} \wk \nu^j_{V \setminus \condv}$ as
  $\level \to \infty $, for any $j$.\ We prove weak convergence by
  proving the pointwise convergence of characteristic functions by
  means of induction on $\{1,\dots,N\}$.\ The claim is immediate for
  $j=1$ from Assumption $\BA$.\ Assume the claim is true for a $j$
  with $1<j<N$.\ Extend the graph by including the clique $C_{j+1}$.

  The relevant characteristic functions are then
  \[
    \chi^{t,j+1}(\bm \omega_{V^{j+1}}) = \smash{
      \int_{\RR^{|V^{j+1}|}} \exp \left (\iota \,\bm
        \omega_{V^{j+1}}.\bm y_{V^{j+1}} \right) } \nu_{C_{j+1} \sm
      S_{j+1}|S_{j}}(d \bm y_{C_{j+1} \sm S_{j+1}}) \nu^{t,j}_{V }(d
    \bm y_{V \sm v}),
  \]
  and 
  \[
    \chi^{j+1}(\bm \omega_{V^{j+1}}) = \smash{ \int_{\RR^{|V^{j+1}|}}
      \exp \left (\iota \,\bm \omega_{V^{j+1}}.\bm y_{V^{j+1}} \right)
    } G^{S_{j+1}}_{C_{j+1} \sm S_{j+1}}(d \bm y_{C_{j+1} \sm S_{j+1}})
    \nu^{j}_{V }(d \bm y_{V \sm v}).
  \]
  We also define a characteristic function associated with the
  conditional probability measure
  \[
    \chi^{\ttime+1 | \ttime, t}(\bm \omega_{C_{\ttime+1} \sm
      S_{\ttime+1}}, y_{S_1}. \dots, y_{S_{\ttime}}) = \smash{
      \int_{\RR^{|C_{\ttime+1} \sm S_{\ttime+1} |}} } \exp \left
      (\iota \,\bm \omega_{C_{\ttime+1} \sm S_{\ttime+1}}.\bm
      y_{C_{\ttime+1} \sm S_{\ttime+1}} \right) \nu_{C_{\ttime} \sm
      S_{\ttime}|S_{\ttime-1}}(d \bm y_{C_{\ttime} \sm S_{\ttime}}).
  \]
  Now,
  \begin{IEEEeqnarray*}{rCl}
    \chi^{t,j+1}(\omega_{V^{j+1}}) &=& \smash{ \int_{\RR^{|V^{j}|}}}  \exp \left (\iota \,\bm  \omega_{V^{j}}.\bm y_{V^{j}} \right) \\ \\
    & \times & \smash{ \int_{\RR^{|V^{|C_{j+1} \sm S_{j+1}|}|}}} \exp
    \left (\iota \,\bm \omega_{C_{j+1} \sm S_{j+1}}.\bm y_{C_{j+1} \sm
        S_{j+1}} \right) \nu_{C_{j+1} \sm S_{j+1}|S_{j}}(d \bm
    y_{C_{j+1} \sm
      S_{j+1}})  \nu^{t,j}_{V }(d \bm y_{V \sm v})\\ \\
    &=& \smash{ \int_{\RR^{|V^{j}|}}} \exp \left (\iota \,\bm
      \omega_{V^{j}}.\bm y_{V^{j}} \right) \chi^{j+1 | j, t}(\bm
    \omega_{C_{j+1} \sm S_{j+1}}, y_{S_1}. \dots, y_{S_{j}})
    \nu^{t,j}_{V }(d \bm y_{V \sm v}).
  \end{IEEEeqnarray*}
  From Assumption $\BA$ and from Lemma \ref{lemma:inf} we have that
  for any $(y_{S_1}, \dots, y_{S_{j}})$, as $\level \to \infty$,
  \begin{equation*}
    \chi^{j+1 | j, t}(\bm \omega_{C_{j+1} \sm S_{j+1}}, y_{S_1}. \dots, y_{S_{j}}) \rightarrow
    \int_{\RR^{|C_{j+1} \sm S_{j+1} |}}  \exp \left (\iota \,\bm  \omega_{C_{j+1} \sm S_{j+1}}.\bm y_{C_{j+1} \sm S_{j+1}} \right) G^{S_{j+1}}_{C_{j+1} \sm S_{j+1}|}(d \bm y_{C_{j+1} \sm S_{j+1}}).
  \end{equation*}
  Moreover, from Lemma \ref{lemma:inf}, this latter convergence is
  uniform on compact sets of
  $(y_{S_2}, \dots y_{S_{j +1}} ) \in \RR^{j}$.\ Hence, from Lemma
  \ref{lemma:ks}, it follows that
  $\chi^{t,j+1}(\omega_{V^{j+1}}) \rightarrow
  \chi^{j+1}(\omega_{V^{j+1}}).$
  \section{Proofs for specific examples}
  \label{app:exproofs}
  \subsection{Decomposable graphical models with
    H\"{u}sler--Reiss cliques}
  \begin{proof}[Proof of Lemma \ref{lemma:hra2}]
    Let $S\neq \emptyset$ with $S\subset C$.\ The density
    $\lambda^{(C)}(\bm y)$ of the exponent measure $\Lambda^{(C)}(\bm y)$ can
    be written, for any $i \in C$, as
    \begin{equation*}
      \lambda^{(C)}(\bm y) = y_i^{-2} \prod_{\substack{j \in C\sm i}} y_j^{-1} \varphi_{|C|-1} \left \{ \log \left (  \frac{\bm y_{C \sm i}}{y_i}\right) + \frac{\bm \Gamma_{C \sm i ,i }}{2} , \bm \Sigma^{(i)}_C  \right \},\quad i\in C,
    \end{equation*}
    where
    $\bm \Sigma_C^{(i)}= \bm \Gamma_{C\sm s, i}^{(C)} \, \bm 1^\top + (\bm
    \Gamma_{C\sm s, i}^{(C)}\, \bm 1^\top)^\top + \bm \Gamma_{C\sm s, C\sm
      i}^{(C)}$.\ In what follows, we assume that $s=v$ if $v\in C$ and
    $s$ equals an arbitrary vertex in $S$, otherwise.

    We have, for $\bm y_C \in \EE^C$,
    \begin{IEEEeqnarray*}{rCl}
      \Lambda_S^{(C)}(\bm y_C) &=& \frac{\partial^{|S|}}{\partial \bm y_S}   \lambda^{(C)}(\bm y) = \int_{[\bm 0,\, \bm y_{C \sm S}]} \lambda_S^{(C)} \left (\bm x_{C \sm S}, \bm y_S \right ) d \bm x_{C \sm S},\\\\
      &=& y_i^{-2} \prod_{j \in S \sm s} y_j^{-1} \frac{\lvert \bm
        Q_C^{(s)} \rvert^{1/2}} {(2 \pi)^{(|C|-1)/2}} \int_{[\bm 0, \bm
        y_{C \sm S}]}
      \exp\left ( -\bm z_C^\top\, \bm Q_C^{(s)} \,\bm z_C
        / 2\right) \, \left\lvert\frac{\partial \bm z_{C\sm S}}{\partial \bm x_{C\sm S}} \right\lvert\,\dd
      \bm x_{C \sm S},
    \end{IEEEeqnarray*}
    where
    $\bm z_C =(\bm z_{S\sm s}, \bm z_{C\sm S}) \equiv (\bm z_{S\sm
      i}(\bm y_{S\sm s}), \bm z_{C\sm S}(\bm x_{C\sm S}))$ with
    \[ (\bm z_{S\sm s}(\bm y_{S\sm s}), \bm z_{C\sm S}(\bm
      x_{C\sm S}))=\left(\log \left ( \bm y_{S \sm s}/y_s\right) +
        \bm \Gamma_{S \sm s , s }\,/\,2, \log \left ( \bm x_{C \sm
            S}/ y_s\right) + \bm \Gamma_{C \sm S , s }\,/\,2 \right),
    \]
    and
    $\bm Q_C^{(s)} = (\bm \Sigma_C^{(s)})^{-1}$.\ Expanding the
    quadratic form into elements indexed by $C\sm S$ and $S\sm s$,
    we get, after applying the change of variables
    $\bm w_{C\sm S} = \bm z_{C\sm S} (\bm x_{C\sm S})$, that
    \begin{IEEEeqnarray*}{rCl}
      \Lambda_S^{(C)}(\bm y_C) &=& y_s^{-2} \prod_{j \in S \sm s}
      y_j^{-1}\frac{\lvert \bm Q_C^{(s)} \rvert^{1/2}} {(2
        \pi)^{(|C|-1)/2}} \, \exp \left ( -\frac{1}{2}\bm z_{S\sm
          i}^\top\,
        \bm Q_{S\sm s}^{(s)} \, \bm z_{S\sm s} \right ) \times \\\\
      &\times&\int_{\left[-\bm \infty,\,\bm z_{C\sm S}(\bm y_{C\sm S})\right
        ]} \exp \left [ - \bm z_{S\sm s}^\top\, \bm Q_{S\sm s, C
          \sm S} \,\bm w_{C \sm S} -\{\bm w_{C \sm S}^\top \, \bm
        Q_{C \sm S} ^{(s)} \, \bm w_{C \sm S }\,/\,2\} \right ] d
      \bm w_{C \sm S}.
    \end{IEEEeqnarray*}
    Completing the square in the argument of the integrand gives after
    simplification that
    \begin{IEEEeqnarray*}{rCl} \Lambda_S^{(C)}(\bm y_C) &=& y_s^{-2}
      \prod_{j \in S \sm s} y_j^{-1} \frac{(2 \pi)^{(|C\sm S|-1)/2}
        \lvert \bm Q_C^{(s)} \rvert^{1/2}} {(2 \pi)^{(|C|-1)/2}\lvert
        \bm Q_{C\sm S}^{(s)} \rvert^{1/2}} \Phi \left \{ \log \left (
          \frac{\bm y_{C \sm S}}{y_s} \right) + \frac{\bm \Gamma_{C
            \sm S}}{2};\, \bm \mu^{(s)}(\bm y_{S\sm s}) , \bm
        \Sigma^{(s)}_{ C \sm S, C \sm S} \right \} \times\nonumber \\\\
      &\times& \exp \left \{ \frac{1} {2} {\bm z_{S \sm s}^\top\, [\bm
          Q^{(s)}_{S\sm s, C \sm S} ( \bm Q^{(s)}_{ C \sm S}
          )^{-1} \bm Q^{(s)}_{C \sm S,S\sm s} - \bm Q_{S\sm
            i,S\sm s}^{(s)}]\,\bm z_{S \sm s} } \right \},
    \end{IEEEeqnarray*}
    where
    $\bm \mu^{(s)}(\bm y_{S\sm s}) = - \big( \bm Q^{(s)}_{ C \sm
      S}\big)^{-1} \bm Q^{(s)}_{ C \sm S, S\sm s}\, \bm z_{S\sm s}(\bm
    y_{S\sm s}).\ $ Therefore, after further simplification and upon
    division by $\Lambda_S^{(S)}(\bm y_S)$, we arrive at
    \begin{IEEEeqnarray*}{rCl}
      \frac {\Lambda_S^{(C)}(\bm y_C)} { \Lambda_S^{(S)}(\bm y_S)} &=& 
      \Phi \left \{ \log \left ( \bm y_{C \sm S} \right) ;\, \ \bm
        a_{C\sm S}^{(S)}(\log(\bm y_S))+ 
        \bm m_{C \setminus S}, (\bm Q^{(s)}_{C \sm S})^{-1} \right \},
    \end{IEEEeqnarray*}
    where $\bm a_{C\sm S}^{(S)}$ and $\bm m_{C \setminus S}$ are
    defined in Section \eqref{eq:HR_cliques_example}.
  \end{proof}
  \begin{proof}[Proof of Proposition \ref{prop:HR}]
    The recurrence relationship for the mean is obtained by taking
    expectations on both sides of the tail graphical model equation
    \eqref{eqn:tgexp}.

    For the precision matrix, we rely on a result concerning the
    inversion of matrices with sparse inverse \citep[Theorem
    1]{JOHNSONlocalinversion}.  Let $\bm X_{V'}$ be an arbitrary
    Gaussian random variable with conditional independences described
    by a chordal graph $\GG'$ with vertices $V'$.  Suppose that the
    cliques and separators of $\GG'$ are denoted by
    $\{ C'_{\ttime} : \ttime= 1, \dots,{N} \}$ and
    $\{S'_{\ttime} : \ttime= 2, \dots,{N} \}$, respectively.  Suppose
    also that the covariance matrices for the random subvectors
    $ \bm X_{C'_\ttime} $ and $ \bm X_{S'_\ttime} $ are given by
    $\bm \Sigma_{C'_{\ttime}}$ and $\bm \Sigma_{S'_{\ttime}}$
    respectively, then the precision matrix of the random vector
    $\bm X_{V'}$ is given by
    \begin{equation}
      \bm Q_{V'} =   \sum_{\ttime=1}^N  \breve { \bm Q}^{\ttime}  -  \sum_{\ttime=2}^N  \breve{\bm P}^{\ttime-1}, 
      \label{eqnlocalinv}
    \end{equation}
    where $ \bm Q^{\ttime} = (\bm \Sigma_{C'_{\ttime}})^{-1}$ and
    ${\bm P}^{\ttime-1} = (\bm \Sigma_{S'_{\ttime}})^{-1}$.\ Next,
    derive the matrices $\bm Q^{\ttime}$ and $\bm P^{\ttime -1}$ for
    the Gaussian tail graphical model
    $\bm Z^{(\condv)}_{ V \sm \condv}$ 
    using an iterative procedure.\ When
    $C_{\ttime}\cap \{v\} \neq \emptyset$, \ we have
    $ \bm Q^{\ttime} = (\bm \Sigma^{(\condv)}_{C_{\ttime} \sm
      \condv})^{-1}$.  In particular, we assume throughout that
    $C_{1}\cap \{\condv\} \neq \emptyset$, so
    $ \bm Q^{1} = (\bm \Sigma^{(\condv)}_{C_{1} \sm \condv})^{-1}$.\
    For $\ttime >1$, the separator $S_{\ttime}$ is a subgraph of both
    $C_{\ttime}$ and $C_{\ttime-1}$.  Hence given $\bm Q^{\ttime-1}$
    we can derive ${\bm P}^{\ttime-1}$, which is the the precision of
    the Gaussian random vector $\bm Z^{(\condv)}_{S_{\ttime}} $, using
    the formula
    \[{\bm P}^{\ttime-1} = \left [ \left ( \left \{ \bm Q^{\ttime-1}
          \right \} ^{-1} \right ) _{S_{\ttime},S_{\ttime} } \right
      ]^{-1}.
    \]

    When $C_{\ttime}\cap \{v\} = \emptyset$, Lemma \ref{lemma:hra2}
    implies that the Gaussian random vector
    $\bigg(\bm Z^{(\condv)}_{C_{\ttime} \sm S_{\ttime} } - \bm
    A_{C_{\ttime} \sm S_{\ttime}, S_{\ttime}}\bm
    Z^{(\condv)}_{S_{\ttime} } , \bm Z^{(\condv)}_{S_{\ttime}} \bigg )
    $ has a distribution with precision
    \begin{equation*}
      \left [
        \begin{array}{cc}
          \bm Q^{(s_{\ttime})}_{C_{\ttime} \sm S_{\ttime}} & 0 \\
          0 & \bm  P^{{\ttime-1}}
        \end{array} \right ].
    \end{equation*} 
    Hence the Gaussian random vector
    $ \bm Z^{(\condv)}_{C_\ttime} = \left (\bm
      Z^{(\condv)}_{C_{\ttime} \sm S_{\ttime} } , \bm
      Z^{(\condv)}_{S_{\ttime}} \right ) $ has a distribution with
    precision
    \begin{equation*}
      \bm Q^{\ttime} = \left [
        \begin{array}{cc}
          \bm Q^{(s_{\ttime})}_{C_{\ttime} \sm S_{\ttime}} +   \bm A_{C_{\ttime} \sm S_{\ttime}, S_{\ttime}} {\bm P}^{\ttime-1} \Big (\bm A_{C_{\ttime}\sm S_{\ttime}, S_{\ttime}} \Big )^T & -  \bm A_{C_{\ttime}\sm S_{\ttime}, S_{\ttime}} {\bm P}^{\ttime-1}\\[8pt]
          -{\bm P}^{\ttime-1}\Big ( \bm A_{C_{\ttime}\sm S_{\ttime}, S_{\ttime}} \Big  )^T &  {\bm P}^{\ttime-1}
        \end{array} \right ],
    \end{equation*}
    and the result for the precision of the distribution of
    $\bm Z^{(\condv)}_{V \sm \condv}$ follows from
    \eqref{eqnlocalinv}.
  \end{proof}
\end{appendices}

\subsection*{Acknowledgments}
A.C and I.P would like to thank Graeme Auld for reading earlier parts
of the manuscript and giving helpful comments that helped improve the
clarity of the paper.
\bibliographystyle{agsm}

\end{document}